\documentclass[11pt,reqno]{amsart}
\usepackage{graphicx}
\usepackage{verbatim}
\usepackage{textcomp}
\usepackage{amssymb}
\usepackage{cite}
\usepackage{amsmath}
\usepackage{latexsym}
\usepackage{amscd}
\usepackage{amsthm}
\usepackage{mathrsfs}
\usepackage{xypic}
\usepackage{bm}
\usepackage{url}
\usepackage{hyperref}

\newtheorem{thm}{Theorem}[section]

\theoremstyle{definition}

\theoremstyle{remark}

\numberwithin{equation}{section} 

\begin{document}
\title[An immersed $S^{n}$ $\lambda$-hypersurface]
{An immersed  $S^{n}$ $\lambda$-hypersurface}
\author{ Zhi Li}
\address{School of Mathematical Sciences, South China Normal University, Guangzhou 510631, P.R. China} \email{lizhihnsd@126.com}
\author{Guoxin Wei}
\address{School of Mathematical Sciences, South China Normal University, Guangzhou 510631, P.R. China} \email{weiguoxin@tsinghua.org.cn}


\thanks{The second author was partly supported by grant No.
11771154 of NSFC and by Guangdong Province Universities and Colleges Pearl River Scholar Funded Scheme (2018).} \maketitle

\begin{abstract}
In this paper, we construct an immersed, non-embedded $S^{n}$ $\lambda$-hypersurface  in Euclidean spaces $\mathbb{R}^{n+1}$.
\end{abstract}

%\MSC{53C44, 53C42}

\section{Introduction}

\noindent
In the study of mean curvature flow, the important class of solutions are those in which the hypersurface evolves under self-similar shrinking. Indeed, under general mean curvature flow, singularities often develop which can be modeled using self-shrinking solutions.  Such solutions can be identified with a single time-slice of the flow, which gives us a hypersurface called a self-shrinker. Let $X : M^{n} \to \mathbb{R}^{n+1}$ be a smooth $n$-dimensional immersed hypersurface in the
$(n+1)$-dimensional Euclidean space $\mathbb{R}^{n+1}$.
An $n$-dimensional hypersurface $M^{n} \subset \mathbb{R}^{n+1}$ is called a self-shrinker if it satisfies the equation
\begin{equation}
\label{1.1}
H + \langle X,\vec{N} \rangle = 0,
\end{equation}
where $H$ is the mean curvature of the hypersurface $M^n$, $X$ is the position vector  and $\vec{N}$ is the unit normal vector of the hypersurface.

\noindent
In \cite{CW0}, Cheng and Wei introduced the notation of $\lambda$-hypersurfaces by studying the weighted volume-preserving mean curvature flow,
which is defined as the following: a family $X(\cdot, t)$ of smooth immersions
$X(\cdot, t) : M^{n} \to \mathbb{R}^{n+1}$
with $X(\cdot, 0) = X(\cdot)$ is called a weighted volume-preserving mean curvature flow if
$$\frac{\partial X(t)}{\partial t} = -\alpha(t)\vec{N}(t) + \vec{H}(t)$$
holds, where
$$\alpha(t)=\frac{\int_{M} H(t)\langle \vec{N}(t), \vec{N}\rangle e^{-\frac{|X|^{2}}{2}}d\mu}{\int_{M}\langle \vec{N}(t), \vec{N}\rangle e^{-\frac{|X|^{2}}{2}}d\mu},$$
$\vec{H}(t) = \vec{H}(\cdot, t)$ and $\vec{N}(t)$ denote the mean curvature vector and the unit normal
vector of hypersurface $M_{t} = X(M^{n}, t)$ at point $X(\cdot, t)$, respectively and $\vec{N}$ is the
unit normal vector of $X : M^{n} \to \mathbb{R}^{n+1}$. One can prove that the flow preserves the weighted volume $V(t)$ defined by
$V(t) = \int_{M}\langle X(t), \vec{N}\rangle e^{-\frac{|X|^{2}}{2}}d\mu$.
The weighted area functional $A: (-\varepsilon, \varepsilon) \to \mathbb{R}$ is defined by
$A(t) =\int_{M}e^{-\frac{|X(t)|^{2}}{2}}d\mu_{t}$,
where $d\mu_{t}$ is the area element of $M$ in the metric induced by $X(t)$. Let $X(t) : M \to \mathbb{R}^{n+1}$ with $X(0) = X$ be a variation of $X$. If $V(t)$ is constant for any $t$, we call $X(t) : M \to \mathbb{R}^{n+1}$ is a weighted volume-preserving variation of $X$. Cheng and Wei \cite{CW0}
have proved that $ X : M \to \mathbb{R}^{n+1}$  is a critical point of the weighted area functional
$A(t)$  for all weighted volume-preserving variations if and only if there exists constant
$\lambda$ such that the hypersurfaces satisfy the equation
\begin{equation}
\label{1.2}
H + \langle X,\vec{N} \rangle = \lambda.
\end{equation}

\noindent
An immersed hypersurface $X(t): M^{n} \to \mathbb{R}^{n+1}$ is called a $\lambda$-hypersurfaces if the equation \ref{1.2} is satisfied.
Moreover, they defined a $F$-functional of $\lambda$-hypersurfaces and studied $F$-stability, which extended a result of Colding-Minicozzi \cite{CM}.

\noindent
From the definition, we know that if $\lambda = 0$, $X : M^{n} \to \mathbb{R}^{n+1}$ is a self-
shrinker of mean curvature flow, that is, $\lambda$-hypersurfaces is the generaliztion of self-shrinkers.
$\lambda$-hypersurfaces were also studied by McGonagle and Ross in \cite{MR}, where they investigate the following isoperimetric type problem in a Gaussian weighted Euclidean space. It turns out that critical points of this variational problem are $\lambda$-hypersurfaces and the only smooth stable ones are hyperplanes.
More information on $\lambda$-hypersurfaces can be found in \cite{CW1}, \cite{G}, \cite{WXZ}and \cite{ZFC}.

\noindent
The simplest examples of self-shrinkers in $\mathbb{R}^{n+1}$ are the sphere of radius $\sqrt{n}$ centered at the origin, cylinders with an axis through the origin and radius $\sqrt{k}$, $1\leq k\leq n-1$, and planes through the origin. In 1992, Angenent~\cite{A} constructed an embedded self-shrinker which is diffeomorphic to $\mathbb{S}^1 \times \mathbb{S}^{n-1}$.
In 1994, Chopp\cite{Ch} described a numerical algorithm to compute surfaces that are approximately self-similar under mean curvature flow.
Kapouleas, Kleene and M{\o}ller\cite{KKM} given the first rigorous construction of complete, embedded
self-shrinking hypersurfaces under mean curvature flow.
Besides, Drugan and Kleene \cite{DK} presented a new family of non-compact properly
embedded, self-shrinking, asymptotically conical, positive mean curvature ends that are hypersurfaces of revolution with circular
boundaries and proved the important classilcation result. In \cite{M}, M{\o}ller proved that there is a rigorous construction of closed, embedded, smooth mean curvature self-shrinkers with high genus, embedded in Euclidean space $\mathbb{R}^3$. Drugan \cite{D} constructed an immersed and non-embedded $S^{2}$ self-shrinker in $\mathbb{R}^3$. In 2015, Cheng and Wei \cite{CW} constructed an embedded $\lambda$-torus.
More recently, Drugan and Nguyen\cite{DN} used variational methods and a modified curvature
flow to give an alternative proof of the existence of a self-shrinking torus under mean curvature flow.
Drugan, Lee and Nguyen\cite{DLN} surveyed known results on closed self-shrinkers for mean curvature
ow and discuss techniques used in recent constructions of closed self-shrinkers with classical rotational symmetry.
Ross\cite{R} showed that there exists a closed, embedded $\lambda$-hypersurface which is diffeomorhic to $\mathbb{S}^{n-1} \times \mathbb{S}^{n-1} \times \mathbb{S}^{1}$ in $\mathbb{R}^{2n}$.

\noindent
In this paper, inspired by Drugan's paper \cite{D}, we construct an immersed, non-embedded $S^{n}$ $\lambda$-hypersurfaces in $\mathbb{R}^{n+1}$.
\begin{thm}\label{thm1-1}
For small $\lambda < 0$, there exists an immersion $\lambda$-hypersurface
$X: S^n\to \mathbb{R}^{n+1}$  which is not embedding.
\end{thm}

\noindent
The basic idea of the proof of Theorem~\ref{thm1-1} is to construct a curve $\Gamma(s)$ to the geodesic
shooting problem with self-intersections whose rotation about the $x$-axis is an topologically $S^{n}$.
In fact, using the symmetry of the curve equation
with respect to refections across the $z$-axis, it is sufficient to find a curve $\Gamma$ that
intersects the $z$-axis perpendicularly.
Using comparison arguments, we give a detailed description of the first two branches of the curve $\Gamma$ when
the initial height is small.
Following the approach of Angenent in~\cite{A}, we use the continuity argument to find an initial condition that corresponds to a solution whose rotation about the $x$-axis is an immersed and non-embedded $S^{n}$ $\lambda$-hypersurface.

\vspace*{2mm}
%%%%%%%%%%%%%%%%%%%%%%%%%%%%%%%%%
\section{Equations of rotational $\lambda$-hypersurfaces in $\mathbb{R}^{n+1}$}\label{sec-notation}

\noindent
Let $\Gamma(s)=(x(s), z(s)),\ s\in (a_{1}, a_{2})$  be a curve with $z > 0$ in the upper half plane
$\mathbb{H} = \{x+iz \mid z > 0, \ x \in \mathbb{R}, \ i = \sqrt{-1}\}$, where $s$ is arc length parameter of $\Gamma(s)$. We
consider a rotational hypersurface $X:(a_{1}, a_{2}) \times S^{n-1}(1) \rightarrow \mathbb{R}^{n+1}$ in $\mathbb{R}^{n+1}$ defined by
$$
X:(a_{1}, a_{2}) \times S^{n-1}(1) \rightarrow \mathbb{R}^{n+1},\quad X(s, \alpha)=(x(s), z(s)\alpha)\in \mathbb{R}^{n+1},
$$
where $S^{n-1}(1)$ is the $(n-1)$-dimensional unit sphere.

\noindent
By a direct calculation, we get the unit normal vector and the mean curvature as follows:
$$\vec{N} = (-\frac{dz(s)}{ds}, \frac{dx(s)}{ds}\alpha),\quad H
= -\frac{d^{2}x(s)}{ds^{2}}\frac{dz(s)}{ds} + \frac{dx(s)}{ds}\frac{d^{2}z(s)}{ds^{2}} - \frac{n-1}{z(s)}\frac{dx(s)}{ds}.$$
And then,
$$\langle X, \vec{N}\rangle=-x(s)\frac{dz(s)}{ds}+z(s)\frac{dx(s)}{ds}.$$
Since $s$ is arc length parameter of the profile curve $\Gamma(s)=(x(s), z(s))$, we have
$(\frac{dx(s)}{ds})^{2} + (\frac{dz(s)}{ds})^{2} = 1$.
Thus, it follows that
$\frac{dx}{ds}\frac{d^{2}x}{ds^{2}} + \frac{dz}{ds}\frac{d^{2}z}{ds^{2}} = 0$.

\noindent
Therefore, it follows that  $X:(a_{1}, a_{2}) \times S^{n-1}(1) \rightarrow \mathbb{R}^{n+1}$ is a $\lambda$-hypersurface in $\mathbb{R}^{n+1}$ if and only if
\begin{equation}\label{2.1}
-\frac{d^{2}x(s)}{ds^{2}}\frac{dz(s)}{ds}+\frac{dx(s)}{ds}\frac{d^{2}z(s)}{ds^{2}}+\left(z(s) - \frac{n-1}{z(s)}\right)\frac{dx(s)}{ds} - x(s)\frac{dz(s)}{ds} =\lambda.
\end{equation}
That is,
\begin{equation}\label{2.2}
\frac{\frac{d^{2}x(s)}{ds^{2}}}{\frac{dz(s)}{ds}} = \left(z(s) - \frac{n-1}{z(s)}\right)\frac{dx(s)}{ds} - x(s)\frac{dz(s)}{ds} - \lambda,
\end{equation}
where $\frac{dx}{ds}\frac{d^{2}x}{ds^{2}} + \frac{dz}{ds}\frac{d^{2}z}{ds^{2}} = 0$.

\noindent
Taking the third derivative of (\ref{2.2}), we have the following equation,
\begin{equation}\label{2.3}\aligned
\frac{\frac{d^{3}x(s)}{ds^{3}}}{\frac{dz(s)}{ds}}
=& - \frac{\frac{dx(s)}{ds}(\frac{d^{2}x(s)}{ds^{2}})^{2}}{(\frac{dz(s)}{ds})^{3}} + \left(z(s) - \frac{n-1}{z(s)}\right)\frac{d^{2}x(s)}{ds^{2}} \\
&+ \frac{n-1}{z(s)^{2}}\frac{dx(s)}{ds}\frac{dz(s)}{ds} + \frac{x(s)\frac{dx(s)}{ds}\frac{d^{2}x(s)}{ds^{2}}}{\frac{dz(s)}{ds}},
\endaligned\end{equation}
where $\frac{dx}{ds}\frac{d^{2}x}{ds^{2}} + \frac{dz}{ds}\frac{d^{2}z}{ds^{2}} = 0$.

\noindent
Obviously, there are several special solutions of (\ref{2.1}):

(1) $(x, z)=(0, s)$ is a solution.
This curve corresponds to the hyperplane through $(0, 0)$ and $\lambda=0$.

(2)$(x, z) = (a\cos\frac{s}{a}, a\sin\frac{s}{a})$ is a solution, where $a=\frac{-\lambda+\sqrt{\lambda^{2}+4n}}{2}$.
This circle $x^{2} + z^{2} = a^{2}$
corresponds to a sphere $S^{n}(a)$ with radius $a$.

(3) $(x, z)=(-s, a)$ is a solution, where $a=\frac{-\lambda+\sqrt{\lambda^{2}+4(n-1)}}{2}$.
This straight line corresponds to a cylinder $S^{n-1}(a)\times \mathbb{R}$.

\vspace*{2mm}
%%%%%%%%%%%%%%%%%%%%%%%%%%%%%%%%%
\section{The first part of the profile curves}\label{sec-notation}

\noindent
In this section, for $b > -(4n+1)\lambda$, we will study the solutions of~({2.2}) with $(x(0),z(0)) = (b, 0)$ and $(\frac{dx(0)}{ds},\frac{dz(0)}{ds}) = (0, 1)$. In~\cite{D}, Drugan have obtained the existence of solutions of self-shrinker near $z=0$ by discussing the existence, uniqueness and continuous dependence of solutions on initial height at $z$ near $0$. That's also true to $\lambda$-hypersurface. After this, we can use some different comparison estimates to describe the basic behavior of the profile curve $\Gamma(s)$. In the end, we complete the section with a detailed description of $\Gamma(s)$ when the initial height $b > -(4n+1)\lambda$ is small.

\subsection{Basic shape of the first branch of the curve $\Gamma(s)$}

For $b > -(4n+1)\lambda$, let $\Gamma(s)$ be the solution of (\ref{2.2}) with $(x(0),z(0)) = (b, 0)$ and $(\frac{dx(0)}{ds},\frac{dz(0)}{ds}) = (0, 1)$. Then
$\frac{d^{2}x(0)}{ds^{2}} = -\frac{1}{n}(b + \lambda)<0$ so that $\Gamma$ starts out concave to the left.
We will describe the basic shape of the curve $\Gamma(s)$ through following several lemmas.

\par\bigskip \noindent{\bf Lemma 3.1.} {\ $\frac{d^{2}x(s)}{ds^{2}} < 0$.}\\
\begin{proof}
Since $\frac{d^{2}x(0)}{ds^{2}} = -\frac{1}{n}(b + \lambda)$, we have that $\frac{d^{2}x(s)}{ds^{2}} < 0$ at $s$ near 0. Assuming $\frac{d^{2}x(s)}{ds^{2}} = 0$ for some $s > 0$. Choosing $\tilde{s}$ so that $\frac{d^{2}x(\tilde{s})}{ds^{2}} = 0$ and $\frac{d^{2}x(s)}{ds^{2}} < 0$ for $s \in [0,\tilde{s})$, then $\frac{d^{3}x(\tilde{s})}{ds^{3}} \geq 0$. Also because $\frac{dx(0)}{ds} = 0$, we have $\frac{dx(\tilde{s})}{ds} < 0$. Using equation (\ref{2.3}), we get $$0 \leq \frac{d^{3}x(\tilde{s})}{ds^{3}} = \frac{n-1}{z(\tilde{s})^{2}}(\frac{dz(\tilde{s})}{ds})^{2}\frac{dx(\tilde{s})}{ds} < 0,$$ it's contradictory. So we have $\frac{d^{2}x(s)}{ds^{2}} < 0$.
\end{proof}

\noindent
Let $s_{1} = s(b) > 0$ be the arc length of the first time at which the unit tangent vector is $(-1, 0)$($\frac{dz(s_{1})}{ds} = 0$).
In~\cite{CW2}, Cheng and Wei prove that an entirely graphic $\lambda$-hypersurface in Euclidean space is a hyperplane,
then $s_{1}$ must exist and $z(s)$ can't take the value at infinity. By Lemma 3.1 and the definition of $s_{1}$, we have that $\frac{dx(s)}{ds} < 0$ and $\frac{dz(s)}{ds} > 0$ in $(0, s_{1})$.
Hence, the curve $\Gamma(s)$ can be
written as a graph of $x = \gamma(z)$, where $0 < z < z (s_{1})$. If $\gamma'(z) =
\frac{dx}{dz} = -\infty$, the
profile curve $\Gamma$ has a horizontal tangent point, in a sense where $\gamma(z)$ blow up. It follows from (\ref{2.2}) that
\begin{equation}\label{2sec1}\aligned
\frac{\gamma''(z)}{1+\gamma'(z)^2} = \left(z - \frac{n-1}{z} \right) \gamma'(z) - \gamma(z) - \lambda\sqrt{1+\gamma'(z)^2},
\endaligned\end{equation}
where $\gamma''(z) = \frac{\frac{d^{2}x(s)}{ds^{2}}}{(\frac{dz(s)}{ds})^{4}}$, $\gamma'(z) = \frac{\frac{dx(s)}{ds}}{\frac{dz(s)}{ds}} < 0$ and $\frac{dz(s)}{ds} = \frac{1}{\sqrt{1+\gamma'(z)^{2}}}>0$.
By Lemma 3.1, we have that $\gamma''(z) < 0$ and $\gamma'(z) \leq 0$ on $[0, z(s_{1}))$. Taking the third and fourth derivatives of (\ref{2sec1}), we have the following equations,
\begin{equation}\label{2sec2}
\frac{\gamma'''}{1+(\gamma')^2}=\frac{2\gamma'(\gamma'')^2}{( 1+(\gamma')^2 )^2}+\left(z
-\frac{n-1}{z} \right)\gamma''+\frac{n-1}{z^2} \gamma'-\frac{\lambda\gamma'\gamma''}{\sqrt{1+(\gamma')^2}}
\end{equation}
and
\begin{equation}\label{2sec3}\aligned
\frac{\gamma^{(4)}}{1+(\gamma')^2}
=&\frac{6\gamma'\gamma''\gamma'''+2(\gamma'')^3}{(1+(\gamma')^2)^2}-\frac{8(\gamma')^2(\gamma'')^3}{(1+(\gamma')^2)^3}
+\left(z-\frac{n-1}{z}\right)\gamma''' - \frac{2(n-1)}{z^3} \gamma'\\
&+\left(1+\frac{2(n-1)}{z^2}\right)\gamma''+\frac{\lambda(\gamma'\gamma'')^{2}}{(1+(\gamma')^2)^{3/2}}-\frac{\lambda((\gamma'')^{2}+\gamma'\gamma''')}{\sqrt{1+(\gamma')^2}}.
\endaligned\end{equation}

\par\bigskip \noindent{\bf Lemma 3.2.} {\ $z(s_{1}) > \frac{-\lambda+\sqrt{\lambda^{2}+4(n-1)}}{2}$.}\\
\begin{proof}
Since $(\frac{dx(s)}{ds})^{2} + (\frac{dz(s)}{ds})^{2} = 1$ and $\frac{dz(s_{1})}{ds} = 0$, we have $\frac{dx(s_{1})}{ds} = -1$.
Supposing $z(s_{1}) = \frac{-\lambda+\sqrt{\lambda^{2}+4(n-1)}}{2}$,
since there exists a special solution of (\ref{2.1}) which $(x,z) = (-s, a)$ with $a = \frac{-\lambda+\sqrt{\lambda^{2}+4(n-1)}}{2}$,
by the existence and uniqueness of solutions for the differential equation, it's contradictory.
Supposing $z(s_{1}) < \frac{-\lambda+\sqrt{\lambda^{2}+4(n-1)}}{2}$, since $\frac{dz(s_{1})}{ds} = 0$, $\frac{dx(s_{1})}{ds} = -1$ and $\frac{dx}{ds}\frac{dx^{2}}{ds^{2}} + \frac{dz}{ds}\frac{dz^{2}}{ds^{2}} = 0$, using equation (\ref{2.1}), we have
\begin{equation}
\frac{d^{2}z(s_{1})}{ds^{2}} = \frac{n-1}{z(s_{1})} - z(s_{1}) - \lambda > 0.\nonumber
\end{equation}
Then, there exists $\delta > 0$ so that $\frac{d^{2}z(s)}{ds^{2}} > 0$ for $s \in (s_{1}-\delta, s_{1})$, and then $\frac{dz(s)}{ds} < 0$ for $s \in (s_{1}-\delta, s_{1})$, it's contradictory. We can get $z(s_{1}) > \frac{-\lambda+\sqrt{\lambda^{2}+4(n-1)}}{2}$.
\end{proof}

\par\bigskip \noindent{\bf Lemma 3.3.} {\ $\lim_{s \to s_{1}}x(s) > -\infty$.}\\
\begin{proof}
When the profile curve $\Gamma$ can be written in the form $(\gamma(z), z)$, $\gamma$ be a solution of (\ref{2sec1}) with $\gamma(0) = b > -(4n+1)\lambda$, $\gamma'(0)=0$.
By Lemma 3.2, we have that $z(s_{1}) > \frac{-\lambda+\sqrt{\lambda^{2}+4(n-1)}}{2}$ and $\lim_{z\to z(s_{1})}\gamma'(z) = -\infty$.
Fixing $0< \delta < 1$ so that $z(s_{1})-\delta > \frac{-\lambda+\sqrt{\lambda^{2}+4(n-1)}}{2}$ and letting $m>0$ so that $z - \frac{n-1}{z} + \lambda \geq m$ for $z \in (z(s_{1}) - \delta, z(s_{1}))$. Choosing $M>0$ so that $m \geq \frac{3}{2M^2}$ and $M \geq -\gamma'(z(s_{1})- \delta)$.

\noindent
For $\delta>\varepsilon>0$, $g_{\varepsilon}(z)$ is defined as
$$g_{\varepsilon}(z) = \frac{M}{\sqrt{(z(s_{1}) - \varepsilon) - z}}, \ \ \ z \in (z(s_{1}) - \delta, z(s_{1}) - \varepsilon).$$
Then we have
$$g_{\varepsilon}(z)> M, \ \ g_{\varepsilon}'(z) = \frac{M}{2((z(s_{1}) - \varepsilon) - z)^{3/2}} > 0$$ and
$$g_{\varepsilon}''(z) = \frac{3}{2M^2} g_{\varepsilon}(z)^2 g_{\varepsilon}'(z) \leq \left(z - \frac{n-1}{z} + \lambda \right) g_{\varepsilon}(z)^2 g_{\varepsilon}'(z),$$ for $z \in (z(s_{1}) - \delta, z(s_{1}) - \varepsilon)$. We use the function $g_{\varepsilon}$ to prove that $-\gamma'$ blow-up no faster than $\frac{M}{\sqrt{z(s_{1})-z}}$.

\noindent
Letting $f(z) = -\gamma'(z)$, then we have $f(z) \geq 0$ and $f'(z) > 0$. Using equation (\ref{2sec2}), we get
\begin{equation}\aligned
f''(z) = & \frac{2f(z)f'(z)^{2}}{1+f(z)^{2}}+\left(z - \frac{n-1}{z} \right)f'(z)(1+f(z)^{2}) \\
& + \frac{(n-1)f(z)}{z^{2}}(1+f(z)^{2}) + \lambda f'(z)f(z)\sqrt{1+f(z)^{2}} \\
\geq & \left(z - \frac{n-1}{z} \right)f'(z)(1+f(z)^{2}) + \lambda f'(z)f(z)\sqrt{1+f(z)^{2}} \nonumber \\
\geq & \left(z - \frac{n-1}{z} + \lambda \right)f'(z)(1+f(z)^{2}) \\
\geq & \left(z - \frac{n-1}{z} + \lambda \right)f'(z)f(z)^{2},
\endaligned\end{equation}
where $\lambda<0, \ \  z \geq \frac{-\lambda+\sqrt{\lambda^{2}+4(n-1)}}{2}$.

\noindent
Next, the purpose is to prove $f \leq g_\varepsilon$. It is known that $$f(z(s_{1})-\delta) \leq M < g_\varepsilon(z(s_{1})-\delta)$$ and $$f(z(s_{1})-\varepsilon) < \lim_{z \to (z(s_{1}) - \varepsilon)} g_\varepsilon(z).$$
Therefore, if there exists some points on $(z(s_{1}) - \delta, z(s_{1}) - \varepsilon)$ so that $f > g_{\varepsilon}$, then $f-g_{\varepsilon}$ achieves a positive maximum at point $\tilde{z} \in (z(s_{1}) - \delta, z(s_{1}) - \varepsilon)$. This leads to $(f-g_{\varepsilon})'(\tilde{z}) = 0$ and $(f-g_{\varepsilon})''(\tilde{z}) \leq 0$. We have $$0 \geq (f-g_{\varepsilon})''(\tilde{z}) \geq \left(\tilde{z} - \frac{n-1}{\tilde{z}} + \lambda \right)f'(\tilde{z}) \left( f(\tilde{z})^2-g_{\varepsilon}(\tilde{z})^2 \right) > 0,$$ it's contradictory.
It follows that $f \leq g_{\varepsilon}$ on $(z(s_{1}) - \delta, z(s_{1}) - \varepsilon)$. Taking $\varepsilon \to 0$, we have that $$\gamma'(z) \geq \frac{-M}{\sqrt{z(s_{1}) - z}},$$ for $x \in (z(s_{1}) - \delta, z(s_{1}))$.
Integrating the inequality from $z(s_{1}) - \delta$ to $z(s_{1})$,
we have $$\lim_{z \to z(s_{1})}\gamma(z) \geq 2M\sqrt{\delta}+\gamma(z(s_{1}) - \delta)>-\infty.$$
Then we get $\lim_{s \to s_{1}} x(s) > -\infty$.
\end{proof}

\subsection{Estimates for small initial height}
From the above lemmas, we have the basic description of the $\Gamma$ curves: For $b > -(4n+1)\lambda$, when the profile curve $\Gamma$ can be written in the form $(\gamma_{b}(z), z)$, $\gamma_{b}$ be a solution of (\ref{2sec1}) with $\gamma_b(0) = b$, $\gamma_b'(0) = 0$. Then $\gamma_b$ is decreasing, concave down and there exists a point $z(s^b_{1}) \in (\frac{-\lambda + \sqrt{\lambda^{2}+4(n-1)}}{2} , \infty)$ so that $\gamma_b$ is defined on $[0, z(s^b_{1}))$ and $\lim_{z \to z(s^b_{1})} \gamma_b'(z) = -\infty$(in a sense $\gamma_b$ blow up at $z(s^b_{1})$). There also exists a point $x(s^b_{1}) \in (-\infty, b)$ so that $\gamma_b(z(s^b_{1})) = x(s^b_{1})$.
Next, we show estimates for $x(s^b_{1})$ and $z(s^b_{1})$ when the initial height $b > -(4n+1)\lambda$ is small.

\par\bigskip \noindent{\bf Proposition 3.4.} {\ For $b > -(4n+1)\lambda$, let $\gamma_b$ denote the solution of (\ref{2sec1}) with $\gamma_b(0) = b$ and $\gamma_b'(0) = 0$. Let $z(s^b_{1})$ denote the point where $\gamma_{b}$ blow-up and $x(s^b_{1}) = \gamma_b(z(s^b_{1}))$.
There exists $\bar{b} > 0$ so that if $b \in (0,\bar{b}]$,
then $$z(s^b_{1}) \geq \sqrt{\ln{\frac{1}{2\sqrt{\pi}(b+\lambda)}}},$$ $$-\frac{4(3n+1)}{\sqrt{\ln{\frac{1}{2\sqrt{\pi}(b+\lambda)}}} + 8\lambda} \leq x(s^b_{1}) < 0,$$ and there exists a point $z_{0}^b \in [\frac{-\lambda+\sqrt{\lambda^{2}+4n}}{2}, \sqrt{\frac{2nb}{b+\lambda}}]$ so that $\gamma_b(z_{0}^b) = 0$.}\\

\noindent
Before we prove this proposition, we will show some lemmas about solutions of (\ref{2sec1}) when the initial height $b > -(4n+1)\lambda$ is small. Let $\gamma$ be the solution of (\ref{2sec1}) with $\gamma(0) = b > -(4n+1)\lambda$ and $\gamma'(0) = 0$.

\par\bigskip \noindent{\bf Lemma 3.5.} {\ Suppose $-(4n+1)\lambda< b <\sqrt{\frac{1}{16\pi e^{36n}}}-\lambda$, then $z(s_{1}) > 3\sqrt{2n}$ and $|\gamma'(z)| \leq \frac{1}{2}$ for $z \in [0, 3\sqrt{2n}]$.}\\
\begin{proof}
Since $\gamma'(0) = 0$, $\gamma'' < 0$ and $\lim_{z\to z(s_{1})} \gamma'(z) = - \infty$, we have that there exists $\tilde{z} \in (0, z(s_{1}))$ so that $\gamma'(\tilde{z}) = -\frac{1}{2}$. For $z \in (0,\tilde{z})$, we get
\begin{equation}\aligned
\frac{d}{dz}\left(e^{-z^2}\gamma'(z)\right)
= &e^{-z^2}\gamma''(z)-2z e^{-z^2}\gamma'(z) \\
\geq &\frac{2}{1 + \gamma'(z)^2} e^{-z^2}\gamma''(z) - 2ze^{-z^2} \gamma'(z)  \\
=& 2e^{-z^2}\left[\left(z-\frac{n-1}{z}\right)\gamma'-\gamma(z)-\lambda\sqrt{1+(\gamma')^{2}}\right] - 2ze^{-z^2}\gamma'(z) \\
=& -2 e^{-z^2}\cdot\frac{n-1}{z}\gamma'(z)-2e^{-z^2}\gamma(z)-2\lambda e^{-z^2}\sqrt{1+(\gamma')^{2}}\nonumber \\
\geq &-e^{-z^2}(2\gamma(z)+2\lambda).
\endaligned\end{equation}

\noindent
Integrating both sides of this inequality from $0$ to $\tilde{z}$,
$$-\frac{1}{2}e^{-(\tilde{z})^2}\geq - \int_0^{\tilde{z}}e^{-z^2}(2\gamma(z)+2\lambda)dz \geq -2(b+\lambda)\int_0^{\tilde{z}} e^{-z^2} dx \geq -2(b+\lambda)\sqrt{\pi},$$
we get $\tilde{z} \geq \sqrt{\ln\frac{1}{4\sqrt{\pi}(b+\lambda)}}$.
When $b < \sqrt{\frac{1}{16\pi e^{36n}}}-\lambda$, we have $e^{-(\tilde{z})^2} < e^{-18n}$, and then $z(s_{1})>\tilde{z} > 3\sqrt{2n}$.
\end{proof}

\par\bigskip \noindent{\bf Lemma 3.6.} {\ If $|\gamma'(z)| \leq \frac{1}{2}$ for $z \in [0, 3\sqrt{2n}]$ and $-(4n+1)\lambda < b < \frac{1}{3\sqrt{2n}}$, then $\frac{z\gamma'(z) - \gamma(z)}{\sqrt{1 + \gamma'(z)^2}}$ is non-increasing on $[0, 3\sqrt{2n}]$.}\\
\begin{proof}
Taking the derivative of $\frac{z\gamma'(z) - \gamma(z)}{\sqrt{1 + \gamma'(z)^2}}$, we have $$ \frac{d}{dz} \left(\frac{z\gamma'(z) - \gamma(z)}{\sqrt{1 + \gamma'(z)^2}} \right) = \gamma''(z) \frac{z + \gamma(z) \gamma'(z)}{(1 + \gamma'(z)^2)^{3/2}}.$$
Next, the purpose is to prove that $z + \gamma(z)\gamma'(z) \geq 0$. Since $z + \gamma(z) \gamma'(z) = 0$ at $z=0$, we only need to prove $1 + \gamma(z) \gamma''(z) + \gamma'(z)^2 \geq 0$ on $(0, 3\sqrt{2n}]$. For $z \in (0, 3\sqrt{2n}]$, according to $|\gamma'(z)| \leq \frac{1}{2}$ and $-(4n+1)\lambda < b < \frac{1}{3\sqrt{2n}}$, we get
\begin{equation}\aligned
\gamma''(z) = & (1 + \gamma'(z)^2) \left[\left(z - \frac{n-1}{z}\right)\gamma'(z) - \gamma(z)
- \lambda\sqrt{1+\gamma'(z)^{2}}\right] \nonumber \\
\geq & (1 + \gamma'(z)^2) \left[z\gamma'(z) - \gamma(z)-\lambda\right] \nonumber \\
\geq & \frac{5}{4} \left[3\sqrt{2n}(-\frac{1}{2}) - b\right] \\
\geq & \frac{5}{4} \left[(-\frac{3}{2})\sqrt{2n} - \frac{1}{3\sqrt{2n}} \right] \geq -3\sqrt{2n}. \nonumber
\endaligned\end{equation}
It follows that $1 + \gamma(z)\gamma''(z) + \gamma'(z)^2 \geq 0$.
\end{proof}

\par\bigskip \noindent{\bf Lemma 3.7.} {\ If $|\gamma'(z)| \leq \frac{1}{2}$ for $z \in [0,3\sqrt{2n}]$ and $-(4n+1)\lambda < b < \frac{1}{3\sqrt{2n}}$, then $\gamma'''(z) < 0$ for $z \in (0, z(s_{1}))$.}\\
\begin{proof}
Using equations (\ref{2sec2}) and (\ref{2sec3}) at $z=0$, we have that $\gamma'''(0) = 0$ and $$\gamma^{(iv)}(0) = \frac{3}{n+2}\gamma''(0)\left(2\gamma''(0)^{2}-\lambda\gamma''(0) + 1\right).$$
And for $\gamma''(0) < 0$ and $-(4n+1)\lambda < b < \frac{1}{3\sqrt{2n}}$, we have that $\gamma^{(iv)}(0) < 0$ and
$\frac{-\lambda+\sqrt{\lambda^{2}+4(n-1)}}{2} < 3\sqrt{2n}$.
Therefore, $\gamma'''(z) < 0$ when $z > 0$ is near 0. Besides, using equation (\ref{2sec2}), we get
\begin{equation}\aligned
\frac{\gamma'''}{1+(\gamma')^2}
= & \frac{2\gamma'(\gamma'')^2}{( 1+(\gamma')^2 )^2}+\left(z - \frac{n-1}{z} \right)\gamma''+\frac{n-1}{z^2} \gamma'-\frac{\lambda\gamma'\gamma''}{\sqrt{1+(\gamma')^2}} \\
< & \frac{2\gamma'(\gamma'')^2}{( 1+(\gamma')^2 )^2}+\left(z
-\frac{n-1}{z} + \lambda\right)\gamma''+\frac{n-1}{z^2} \gamma',\nonumber
\endaligned\end{equation}
where  $0 <\frac{-\gamma'}{\sqrt{1+(\gamma')^2}} < 1$.
It follows that $\gamma'''(z) < 0$ when $z \geq \frac{-\lambda+\sqrt{\lambda^{2}+4(n-1)}}{2}$.
If there are some points on $(0, z(s_{1}))$ so that $\gamma'''(z) = 0$, there must exist $\tilde{z} \in (0,\frac{-\lambda+\sqrt{\lambda^{2}+4(n-1)}}{2})$ so that $\gamma'''(\tilde{z}) = 0$ and $\gamma'''(z) < 0$ for $z \in (0,\tilde{z})$. It follows that $\gamma^{(iv)}(\tilde{z}) \geq 0$. It is obvious that $z \gamma''(z) -\gamma'(z)$ is decreasing and negative on $(0,\tilde{z})$.
By Lemma 3.6, we have $\gamma''(\tilde{z}) \geq -3\sqrt{2n}$.
Then, using equation (\ref{2sec3}) and $|\gamma'(\tilde{z})| \leq \frac{1}{2}$,
we get
\begin{equation}\aligned
\frac{\gamma^{(iv)}(\tilde{z})}{1 + \gamma'(\tilde{z})^2}
= & 2(\gamma''(\tilde{z}))^3\frac{1 - 3\gamma'(\tilde{z})^2}{(1 + \gamma'(\tilde{z})^2)^3} + \gamma''(\tilde{z}) + 2(n-1)\frac{\tilde{z}\gamma''(\tilde{z}) - \gamma'(\tilde{z})}{(\tilde{z})^3} \\
& + \frac{\lambda(\gamma'(\tilde{z})\gamma''(\tilde{z}))^{2}}{(1 + \gamma'(\tilde{z})^2)^{3/2}}
- \frac{\lambda(\gamma''(\tilde{z}))^{2}}{\sqrt{1 + \gamma'(\tilde{z})^2}}\\
= & 2(\gamma''(\tilde{z}))^3\frac{1 - 3\gamma'(\tilde{z})^2}{(1 + \gamma'(\tilde{z})^2)^3} + \gamma''(\tilde{z}) + 2(n-1)\frac{\tilde{z}\gamma''(\tilde{z}) - \gamma'(\tilde{z})}{(\tilde{z})^3} \\
& -\frac{\lambda(\gamma''(\tilde{z}))^{2}}{(1 + \gamma'(\tilde{z})^2)^{3/2}} \\
\leq & 2(\gamma''(\tilde{z}))^3 \frac{1 - 3\gamma'(\tilde{z})^2}{(1+\gamma'(\tilde{z})^2)^3} + \frac{\gamma''(\tilde{z})(1- \lambda\gamma''(\tilde{z}))}{(1 + \gamma'(\tilde{z})^2)^{3/2}} + 2(n-1)\frac{\tilde{z}\gamma''(\tilde{z}) - \gamma'(\tilde{z})}{(\tilde{z})^3}\\
\leq & 2(\gamma''(\tilde{z}))^3 \frac{1 - 3\gamma'(\tilde{z})^2}{(1+\gamma'(\tilde{z})^2)^3} + \frac{\gamma''(\tilde{z})(1+3\sqrt{2n}\lambda)}{(1 + \gamma'(\tilde{z})^2)^{3/2}} + 2(n-1)\frac{\tilde{z}\gamma''(\tilde{z}) - \gamma'(\tilde{z})}{(\tilde{z})^3}\\
< & 0. \nonumber
\endaligned\end{equation}
where $\gamma''(\tilde{z}) < 0$, $|\gamma'(\tilde{z})| \leq \frac{1}{2}$ and $-(4n+1)\lambda < b < \frac{1}{3\sqrt{2n}}$,
it's contradictory. So we have $\gamma'''(z) < 0$ for $z \in (0, z(s_{1}))$.
\end{proof}

\par\bigskip \noindent{\bf Lemma 3.8.} {\ Suppose $-(4n+1)\lambda< b <\sqrt{\frac{1}{16\pi e^{36n}}}-\lambda$, then there exists a point $z_0 \in [\frac{-\lambda+\sqrt{\lambda^{2}+4n}}{2}, \sqrt{\frac{2nb}{b+\lambda}}]$ so that $\gamma(z_0) = 0$.}\\
\begin{proof}
We assume $-(4n+1)\lambda< b <\sqrt{\frac{1}{16\pi e^{36n}}}-\lambda$($b < \frac{1}{3\sqrt{2n}}$), by Lemma 3.5,
we have that $z(s_{1}) > 3\sqrt{2n}$ and $|\gamma'(z)| \leq \frac{1}{2}$ for $x \in [0, 3\sqrt{2n}]$. Therefore, by Lemma 3.7, we have $\gamma'''<0$ on $(0, z(s_{1}))$. Integrating this inequality from $0$ to $z$ repeatedly, we have
$$\gamma(z) < b- \frac{1}{2n}(b+\lambda)z^2.$$

\noindent
If $b- \frac{1}{2n}(b+\lambda)z^2 \leq 0$, we have $\gamma(z) < 0$ for $z \geq \sqrt{\frac{2nb}{b+\lambda}}$.
Since $b > -(4n+1)\lambda$, we have $\sqrt{\frac{2nb}{b+\lambda}} < 2\sqrt{n}$.
Besides, $z(s_{1}) > 3\sqrt{2n}$, there exists $z_{0} \in (0,\sqrt{\frac{2nb}{b+\lambda}})$ so that $\gamma(z_{0}) = 0$.

\noindent
Next, we will evaluate the lower bound of $z_0$. We write equation (\ref{2sec1}) as the form
\begin{equation}\label{2sec4}
\frac{d}{dz} \left( \frac{z^{n-1} \gamma'(z)}{ \sqrt{1+\gamma'(z)^2} } \right)
= z^{n-1}\cdot\frac{z\gamma'(z) - \gamma(z)}{\sqrt{1+\gamma'(z)^2} } - \lambda z^{n-1}.
\end{equation}
From Lemma 3.6, we have that $\frac{z\gamma'(z) - \gamma(z)}{ \sqrt{1+\gamma'(z)^2} } \geq \frac{z_{0} \gamma'(z_{0}) }{\sqrt{1+\gamma'(z_{0})^2} }$ for $z\in [0, z_{0}]$. Integrating (\ref{2sec4}) from $0$ to $z_{0}$, we get
\begin{equation}\aligned
\frac{z_{0}^{n-1} \gamma'(z_{0})}{ \sqrt{1+\gamma'(z_{0})^2} } = & \int_0^{z_{0}} z^{n-1} \cdot \frac{z \gamma'(z) - \gamma(z)}{ \sqrt{1+\gamma'(z)^2} }dz - \int_0^{z_{0}} \lambda z^{n-1}dz  \\
\geq & \frac{z_{0} \gamma'(z_{0})}{ \sqrt{1+\gamma'(z_{0})^2} } \int_0^{z_{0}}z^{n-1}dz - \frac{\lambda}{n}(z_{0})^{n} \\
= & \frac{1}{n}(z_{0})^{n}\frac{z_{0} \gamma'(z_{0})}{ \sqrt{1+\gamma'(z_{0})^2} } - \frac{\lambda}{n}(z_{0})^{n} \\
\geq & \frac{1}{n}(z_{0})^{n}\frac{z_{0} \gamma'(z_{0})}{ \sqrt{1+\gamma'(z_{0})^2}}
+ \frac{\lambda}{n}(z_{0})^{n} \cdot \frac{\gamma'(z_{0})}{ \sqrt{1+\gamma'(z_{0})^2} }. \nonumber
\endaligned\end{equation}
Therefore, we have that $1 \leq \frac{(z_{0})^2}{n} + \frac{\lambda}{n}z_{0}$, then $z_{0} \geq \frac{-\lambda+\sqrt{\lambda^{2}+4n}}{2}$. This proves the last statement of the Lemma.
\end{proof}

\noindent
According to Lemma 2.8,
if $b- \frac{1}{2n}(b+\lambda)z^2 < \lambda$, we have $\gamma(z) < \lambda$ for $z \geq \sqrt{\frac{2n(b-\lambda)}{b+\lambda}}$.
If $b- \frac{1}{2n}(b+\lambda)z^2 < 8n\lambda$, we have $\gamma(z) < 8n\lambda$ for $z \geq \sqrt{\frac{2n(b-8n\lambda)}{b+\lambda}}$.
At this time, $\sqrt{\frac{2n(b-\lambda)}{b+\lambda}} < 2\sqrt{n}$ and $\sqrt{\frac{2n(b-8n\lambda)}{b+\lambda}} < 4\sqrt{n}$ for $b > -(4n+1)\lambda$.

\par\bigskip \noindent{\bf Lemma 2.9.} {\ Suppose $-(4n+1)\lambda < b < \sqrt{\frac{1}{16\pi e^{36n}}}-\lambda$, $z(s_{1})> 3\sqrt{2n}$ and there exists a point $z_{0}\in [\frac{-\lambda+\sqrt{\lambda^{2}+4n}}{2}, \sqrt{\frac{2nb}{b+\lambda}}]$ so that $\gamma(z_{0}) = 0$. Then, for $z \in [z_{0}, z(s_{1}))$, $$\gamma(z) > \frac{12n}{z}\gamma'(z).$$}\\

\begin{proof}
Letting $\Phi(z) = \frac{1}{12n}z\gamma(z) - \gamma'(z)$. To prove $\gamma(z) > \frac{12n}{z}\gamma'(z)$ for $z \in [z_{0}, z(s_{1}))$,
it only need to be satisfied $\Phi(z) > 0$. Since $\Phi(z_0) = -\gamma'(z_0) > 0$ and
\begin{equation}\aligned
\frac{1}{12n}z\gamma(z) = &\frac{1}{12n}z\int_{z_{0}}^z\gamma'(\xi) d\xi \nonumber \\
> &\frac{1}{12n}z(z-z_{0})\gamma'(z), \nonumber
\endaligned\end{equation}
we have $\Phi(z) > \frac{1}{12n}\gamma'(z)(z^{2}-z_{0}z-12n)$.
Besides, for $z_{0} \geq \frac{-\lambda+\sqrt{\lambda^{2}+4n}}{2}$, we have that $\Phi(z) > 0$ when $z\leq \frac{\frac{-\lambda+\sqrt{\lambda^{2}+4n}}{2}+\sqrt{(\frac{-\lambda+\sqrt{\lambda^{2}+4n}}{2})^{2}+48n}}{2}$.

\noindent
Letting $\mathbf{z} = \frac{\frac{-\lambda+\sqrt{\lambda^{2}+4n}}{2}+\sqrt{(\frac{-\lambda+\sqrt{\lambda^{2}+4n}}{2})^{2}+48n}}{2}$,
it is obvious that $\mathbf{z} > 4\sqrt{n}$.
If there are some points $z \in [z_{0}, z(s_{1}))$ so that $\Phi(z) = 0$, then there exists a point $\tilde{z} \in (\mathbf{z}, z(s_{1}))$ so that $\Phi(\tilde{z}) = 0$ and $\Phi(z) > 0$ for $z \in [z_{0}, \tilde{z})$. This shows that $\Phi'(\tilde{z}) \leq 0$ and $\frac{1}{12n}\tilde{z} \gamma(\tilde{z}) = \gamma'(\tilde{z})$.
Since $\gamma(\tilde{z}) < 0$ and $\gamma'(\tilde{z}) < 0$,
we get
\begin{equation}\aligned
0\geq&\Phi'(\tilde{z}) \\
= & \frac{1}{12n} \gamma(\tilde{z}) + \frac{1}{12n}\tilde{z} \gamma'(\tilde{z}) - \gamma''(\tilde{z})  \\
 \geq & \frac{1}{12n} \gamma(\tilde{z}) + \frac{1}{12n}\tilde{z} \gamma'(\tilde{z}) - \frac{\gamma''(\tilde{z})}{1+\gamma'(\tilde{z})^2}  \\
= & \frac{1}{12n}\gamma(\tilde{z}) + \frac{1}{(12n)^{2}}(\tilde{z})^{2}\gamma(\tilde{z}) -
\left[\frac{1}{12n}\tilde{z}\gamma(\tilde{z})\left(\tilde{z}- \frac{n-1}{\tilde{z}}\right) - \gamma(\tilde{z}) - \lambda\sqrt{1+\gamma'(\tilde{z})^2} \right] \\
\geq & \frac{13}{12}\gamma(\tilde{z}) - \frac{12n-1}{(12n)^{2}}(\tilde{z})^2\gamma(\tilde{z}) + \lambda(1 - \gamma'(\tilde{z})) \\
= & \gamma(\tilde{z})\left(\frac{13}{12} -\frac{\lambda}{12n}\tilde{z}- \frac{12n-1}{(12n)^{2}}(\tilde{z})^2 \right) +\lambda. \nonumber
\endaligned\end{equation}

\noindent
Letting $f(z) = \frac{13}{12} - \frac{\lambda}{12n}z- \frac{12n-1}{(12n)^{2}}z^2$ for $z \geq 4\sqrt{n}$.
Since $-(4n+1)\lambda < b < \sqrt{\frac{1}{16\pi e^{36n}}}-\lambda$, we have $-\lambda < \frac{1}{4n}\sqrt{\frac{1}{16\pi e^{36n}}}$. Through the monotonicity of the quadratic function, we have $f(z) < 0$ for $z\geq 4\sqrt{n}$.
Since $\tilde{z} > 4\sqrt{n}$, we have $\gamma(\tilde{z}) < 8n\lambda < 0$, and then
\begin{equation}\aligned
& \gamma(\tilde{z})\left(\frac{13}{12} -\frac{\lambda}{12n}\tilde{z}- \frac{12n-1}{(12n)^{2}}(\tilde{z})^2 \right) +\lambda  \\
>& \lambda\left((\frac{26n}{3}+1) -\frac{2\lambda}{3}\tilde{z}- \frac{12n-1}{18n}(\tilde{z})^2 \right) \\
>& \lambda\left((\frac{26n}{3}+1) -\frac{2\lambda}{3}(4\sqrt{n})- \frac{12n-1}{18n}(4\sqrt{n})^2 \right) \\
=& \lambda\left(-2n + \frac{17}{9} - \frac{8\sqrt{n}\lambda}{3} \right)
> 0, \nonumber
\endaligned\end{equation}
where $-\lambda < \frac{1}{4n}\sqrt{\frac{1}{16\pi e^{36n}}}$ and $n \geq 2$. It's contradictory. So we have that $\gamma(z) > \frac{12n}{z}\gamma'(z)$ for $z \in [z_{0}, z(s_{1}))$.
\end{proof}

\begin{proof}[Proof of Proposition 3.4.]
Let $\gamma$ be the solution of (\ref{2sec1}) with $\gamma(0) = b > 0$ and $\gamma'(0) = 0$. By Lemma 3.2 and Lemma 3.3, there exists a point $z(s_{1}) \in (\frac{-\lambda + \sqrt{\lambda^{2}+4(n-1)}}{2}, \infty)$ and a point $x(s_{1}) \in (-\infty, b)$ so that $\lim_{z \to z(s_{1})} \gamma'(z) = - \infty$ and $\gamma(z(s_{1})) = x(s_{1})$.

\noindent
We want to refine the estimate from Lemma 3.5 to establish a lower bound for $z(s_{1})$ in terms of $b$. Let $z_{1} \in (0, z(s_{1}))$ be the point where $\gamma'(z_{1}) = -1$. Using the same method we used in the proof of Lemma 3.5, integrating both sides of this inequality
$$\frac{d}{dz} \left( e^{-z^2} \gamma'(z) \right) \geq - 2(\gamma(z)+\lambda)e^{-z^2}$$
from $0$ to $z_{1}$. We have $-e^{-(z_{1})^2} \geq -2(b+\lambda)\sqrt{\pi}$
, then $$z_{1} \geq \sqrt{\ln{\frac{1}{2\sqrt{\pi}(b+\lambda)}}}.$$
Assuming $-(4n+1)\lambda\leq b \leq \sqrt{\frac{1}{4\pi e^{64n}}} - \lambda$,
we have $z_{1} \geq 4\sqrt{2n}>\sqrt{\frac{2nb}{b+\lambda}}$, and then $\frac{13n-1}{z} < \frac{3z}{4}$ for $z > z_{1}$.
By Lemma 3.8, there exists $z_{0}\in [\frac{-\lambda+\sqrt{\lambda^{2}+4n}}{2}, \sqrt{\frac{2nb}{b+\lambda}}]$ so that $\gamma(z_{0}) = 0$, then
we have from Lemma 3.9 that $\gamma(z) > \frac{12n}{z}\gamma'(z)$ for $z \in [z_{0}, z(s_{1}))$.
In particular, at $z_{1}$, we get $$\gamma(z_{1}) > - \frac{12n}{z_{1}}.$$

To extend this estimate of $\gamma(z_1)$ to $\gamma(z(s_{1})) = x(s_{1})$. For $z \geq z_{1}$, we get
\begin{equation}\aligned
\gamma''(z) \leq & \gamma'(z)^2 \frac{\gamma''(z)}{1 + \gamma'(z)^2}  \\
= & \gamma'(z)^2 \left[ \left(z - \frac{n-1}{z}\right) \gamma'(z) - \gamma(z) - \lambda\sqrt{1 + \gamma'(z)^2}\right]  \\
< & (z - \frac{13n-1}{z} + 2\lambda)\gamma'(x)^3 \\
\leq & (\frac{1}{4}z + 2\lambda)\gamma'(z)^3, \nonumber
\endaligned\end{equation}
where we have used that $\gamma'(z) \leq -1$ and $\gamma(z) > \frac{12n}{z}\gamma'(z)$ for $z \geq z_{1}$.
Integrating both sides of this inequality from $z$ to $z(s_{1})$,
we get $$\gamma'(z)^2 \leq \frac{4}{(z(s_{1}) - z)(z(s_{1}) + z + 16\lambda)},$$
for $z \geq z_{1}$.
Since $\gamma'(z) < 0$ and $-(4n+1)\lambda < b \leq \sqrt{\frac{1}{4\pi e^{64n}}} - \lambda$, we have that
$$-\frac{1}{4n}\sqrt{\frac{1}{4\pi e^{64n}}}< \lambda <0, \ \ z_{1}+8\lambda>0.$$
Then
\begin{equation}\aligned\label{2sec5}
\gamma'(z) & \geq - \frac{2}{\sqrt{z(s_{1}) - z}\sqrt{z(s_{1}) + z + 16\lambda}} \\
           & \geq - \frac{1}{\sqrt{z(s_{1}) - z}} \cdot \frac{2}{\sqrt{z(s_{1}) + z_{1} + 16\lambda}},
\endaligned\end{equation}
for $z \in [z_{1}, z(s_{1}))$. Taking value at $z_{1}$, we get $$- \frac{\sqrt{z(s_{1}) - z_{1}}}{\sqrt{z(s_{1}) + z_{1} + 16\lambda}} \geq - \frac{2}{z(s_{1}) + z_{1} + 16\lambda}.$$
Integrating (\ref{2sec5}) from $z_{1}$ to $z(s_{1})$, we have $$\gamma(z(s_{1})) - \gamma(z_{1}) \geq -\frac{4}{\sqrt{z(s_{1}) + z_{1} + 16\lambda}} \cdot \sqrt{z(s_{1}) - z_{1}},$$
then
\begin{equation}\aligned
\gamma(z(s_{1})) \geq & \gamma(z_{1}) - \frac{4}{\sqrt{z(s_{1}) + z_{1} + 16\lambda}} \cdot \sqrt{z(s_{1}) - z_{1}}  \\
\geq & - \frac{12n}{z_{1}} -  \frac{8}{z(s_{1}) + z_{1} + 16\lambda}  \\
\geq & -\frac{12n}{z_{1} + 8\lambda} - \frac{8}{z(s_{1}) + z_{1} + 16\lambda}\\
\geq & -\frac{4(3n+1)}{z_{1} + 8\lambda}\\
\geq & -\frac{4(3n+1)}{\sqrt{\ln{\frac{1}{2\sqrt{\pi}(b+\lambda)}}} + 8\lambda}. \nonumber
\endaligned\end{equation}
Finally, we complete the proof of the Proposition 3.4 when $\bar{b} = \sqrt{\frac{1}{4\pi e^{64n}}} - \lambda$.
\end{proof}

\vspace*{2mm}
%%%%%%%%%%%%%%%%%%%%%%%%%%%%%%%%%
\section{The second part of the profile curves}

\noindent
The basic shape of the curve $\Gamma(s)$ is described in the previous section, we have that $\frac{dx(s)}{ds} < 0$
and $\frac{d^{2}x(s)}{ds^{2}} < 0$ when $s\in (0, s_{1})$.
Then the first branch of the curve $\Gamma(s)$ can be written as the curve $(\gamma(z), z)$ where $\gamma$ is the solution to (\ref{2sec1})
with $\gamma(0) = b > -(4n+1)\lambda$ and $\gamma'(0) = 0$.
At $s = s_{1}$, we have $\frac{dx(s_{1})}{ds} = -1$ and $\frac{dz(s_{1})}{ds} = 0$.
According to the equation (\ref{2.2}) and $\frac{dx}{ds}\frac{d^{2}x}{ds^{2}} + \frac{dz}{ds}\frac{d^{2}z}{ds^{2}} = 0$,
we have $\frac{d^{2}z(s_{1})}{ds^{2}} = \frac{n-1}{z(s_{1})} - z(s_{1}) -\lambda < 0$, where $z(s_{1}) > \frac{-\lambda+\sqrt{\lambda^{2}+4(n-1)}}{2}$.
This shows that the curve $\Gamma$ is concave down at $(x(s_{1}), z(s_{1}))$ and heads back towards the $x$-axis as $s$ increases.
Let $s_{2} = s(b) > 0$ be the arc length of the second time, if any, at which either $\frac{dz(s_{2})}{ds} = 0$ or $z(s_{2})= 0$.

\noindent
In this section, we study the curves $\Gamma(s)$ as they travel from $(x(s_{1}), z(s_{1}))$ toward the $x$-axis.
From this definition of $s_{2}$, we have $\frac{dz(s)}{ds} < 0$ for $s \in (s_{1}, s_{2})$.
this curve  can be written as a graph of $x = \beta(z)$, where $z(s_{2}) < z < z(s_{1})$.
It follows from (\ref{2.2}) that
\begin{equation}\label{3sec1}\aligned
\frac{\beta''(z)}{1+\beta'(z)^2} = \left(z - \frac{n-1}{z} \right) \beta'(z) - \beta(z) + \lambda\sqrt{1+\beta'(z)^2},
\endaligned\end{equation}
where $\beta''(z) = \frac{\frac{d^{2}x(s)}{ds^{2}}}{(\frac{dz(s)}{ds})^{4}}$ and $\frac{dz(s)}{ds} = \frac{-1}{\sqrt{1+\beta'(z)^{2}}}< 0$.
Then, $\beta$ is the unique solution to (\ref{3sec1}) with $\beta(z(s_{1})) = x(s_{1})$, $\beta(z(s_{2})) = x(s_{2})$ and $\lim _{z\to z(s_{1})} \beta'(z)= \infty$.
Taking the third derivative of (\ref{3sec1}), we have the following equation,
\begin{equation}\label{3sec2}
\frac{\beta'''}{1+(\beta')^2}=\frac{2\beta'(\beta'')^2}{( 1+(\beta')^2 )^2}+\left(z
-\frac{n-1}{z} \right)\beta''+\frac{n-1}{z^2}\beta'+\frac{\lambda\beta'\beta''}{\sqrt{1+(\beta')^2}}
\end{equation}

\noindent
Since $\frac{dx(s)}{ds}=-1$ near $s=s_{1}$, the curve near the point $(x(s_{1}), z(s_{1}))$ can be written as a graph of $z=\alpha(x)$.
It follows from (\ref{2.1}) that
\begin{equation}\label{3sec3}\aligned
\frac{\alpha''(x)}{1+\alpha'(x)^2} = \left(\frac{n-1}{\alpha} - \alpha\right) + x\alpha'(x) - \lambda\sqrt{1+\alpha'(x)^2},
\endaligned\end{equation}
where $\alpha''(x) = \frac{\frac{d^{2}z(s)}{ds^{2}}}{(\frac{dx(s)}{ds})^{4}}$ and $\frac{dx(s)}{ds} = \frac{-1}{\sqrt{1+\alpha'(x)^{2}}}$.
Taking the third derivative of (\ref{3sec3}), we have the following equation,
\begin{equation}\label{3sec4}
\frac{\alpha'''}{1+(\alpha')^2}=\frac{2\alpha'(\alpha'')^2}{( 1+(\alpha')^2 )^2}
-\frac{n-1}{\alpha^{2}}\alpha'+x\alpha''-\frac{\lambda\alpha'\alpha''}{\sqrt{1+(\alpha')^2}}.
\end{equation}
Next, we will show that for small $b > -(4n+1)\lambda$, $x(s)$ achieves a negative minimum at a point $s_{m} \in (s_{1},s_{2})$, the curve $\Gamma(s)$ is concave to the right, $z(s_{2})>0$ and $0 < x(s_{2}) < \infty$.

\subsection{Basic shape of the second branch of the curve $\Gamma(s)$}
In the first, we give some properties of the second branch of the curve $\Gamma(s)$ that are consequences of equations (\ref{2.2}) and (\ref{3sec1}).

\par\bigskip \noindent{\bf Lemma 4.1.} {\ There exists at most one point $s_{m} \in (s_{1},s_{2})$ so that $\frac{dx(s_{m})}{ds} = 0$. If the point $s_{m}$ exists, then $\frac{d^{2}x(s)}{ds^{2}} > 0$ on $(s_{1}, s_{2})$ and $x(s_{m}) < \lambda$.}\\
\begin{proof}
According to the equation (\ref{2.2}) and $\frac{dx}{ds}\frac{d^{2}x}{ds^{2}} + \frac{dz}{ds}\frac{d^{2}z}{ds^{2}} = 0$,
we have $\frac{d^{2}z(s_{1})}{ds^{2}} = \frac{n-1}{z(s_{1})} - z(s_{1}) -\lambda < 0$,
where $z(s_{1}) > \frac{-\lambda+\sqrt{\lambda^{2}+4(n-1)}}{2}$, $\frac{dx(s_{1})}{ds} = -1$ and $\frac{dz(s_{1})}{ds} = 0$.
There exists $\delta >0$ so that for $s \in (s_{1}, s_{1} + \delta)$, $\frac{dx(s)}{ds} < 0$, $\frac{dz(s)}{ds} < 0$ and $\frac{d^{2}z(s)}{ds^{2}} < 0$.
Using the equation $\frac{dx}{ds}\frac{d^{2}x}{ds^{2}} + \frac{dz}{ds}\frac{d^{2}z}{ds^{2}} = 0$, then
we have $\frac{d^{2}x(s)}{ds^{2}} > 0$ for $s \in (s_{1}, s_{1} + \delta)$.

\noindent
Next, we will prove that there exists at most one point $s_{m} \in (s_{1},s_{2})$ so that $\frac{dx(s_{m})}{ds} = 0$.
Supposing that there are two adjacent zeros $s_{m}$ and $s_{n}$($s_{m} < s_{n}$) to $\frac{dx(s)}{ds}$,
then there exists $\tilde{s}\in (s_{m}, s_{n})$ so that $\frac{d^{2}x(\tilde{s})}{ds^{2}} = 0$. We will discuss this in three cases.
In the first case, if $\frac{dx(s)}{ds} \equiv 0$ for $s \in (s_{m}, s_{n})$, we have $\frac{d^{2}x(s)}{ds^{2}} = 0$ and $\frac{dz(s)}{ds} = -1$ for $s \in (s_{m}, s_{n})$, and then
$x(s) = \lambda$ for $s \in (s_{m}, s_{n})$.
By the uniqueness of solutions for the differential equations, $x$ must be the constant function $x(s) \equiv \lambda$ on $(s_{1}, s_{2})$,
it's contradictory.
In the second case, if $\frac{d^{2}x(s)}{ds^{2}} > 0$ on $(s_{m} , \tilde{s})$ and $\frac{d^{2}x(s)}{ds^{2}} < 0$ on $(\tilde{s}, s_{n})$,
then we have $\frac{dx(\tilde{s})}{ds} > 0$ and $\frac{d^{3}x(\tilde{s})}{ds^{3}} \leq 0$.
Using equation (\ref{2.3}), we get
$$0 \geq \frac{d^{3}x(\tilde{s})}{ds^{3}} = \frac{n-1}{z(\tilde{s})^{2}}(\frac{dz(\tilde{s})}{ds})^{2}\frac{dx(\tilde{s})}{ds} > 0,$$ it's contradictory.
In the third case, if $\frac{d^{2}x(s)}{ds^{2}} < 0$ on $(s_{m} , \tilde{s})$ and $\frac{d^{2}x(s)}{ds^{2}} > 0$ on $(\tilde{s}, s_{n})$,
then we have $\frac{dx(\tilde{s})}{ds} < 0$ and $\frac{d^{3}x(\tilde{s})}{ds^{3}} \geq 0$.
Using equation (\ref{2.3}) again, we get
$$0 \leq \frac{d^{3}x(\tilde{s})}{ds^{3}} = \frac{n-1}{z(\tilde{s})^{2}}(\frac{dz(\tilde{s})}{ds})^{2}\frac{dx(\tilde{s})}{ds} < 0,$$ it's contradictory.

\noindent
Supposing $\frac{dx(s_{m})}{ds} = 0$ for $s_{m} \in (s_{1}, s_{2})$, we have $\frac{dz(s_{m})}{ds} = -1$.
Since $\frac{d^{2}x(s)}{ds^{2}} > 0$ for $s$ close to $s_{1}$ and $\frac{dx(s)}{ds}$ has just one zero point,
it follows that $x(s)$ achieves minimum value at $s_{m}$ and $\frac{d^{2}x(s_{m})}{ds^{2}} > 0$.
Using equation (\ref{2.2}), we have $x(s_{m}) < \lambda$.
Arguing as in the same proof of lemma 3.1,
supposing $\frac{d^{2}x(s)}{ds^{2}}=0$ for some points $s \in (s_{1}, s_{m})$.
Choosing $\bar{s} \in (s_{1}, s_{m})$ so that $\frac{d^{2}x(\bar{s})}{ds^{2}} = 0$ and $\frac{d^{2}x(s)}{ds^{2}} > 0$ for $s \in (\bar{s}, s_{m})$, then we have that $\frac{dx(\bar{s})}{ds} < 0$ and $\frac{d^{3}x(\bar{s})}{ds^{3}} \geq 0$.
Using equation (\ref{2.3}),
we get $$0 \leq \frac{d^{3}x(\bar{s})}{ds^{3}} = \frac{n-1}{z(\bar{s})^{2}}(\frac{dz(\bar{s})}{ds})^{2}\frac{dx(\bar{s})}{ds} < 0,$$ it's contradictory.
The same method, supposing $\frac{d^{2}x(s)}{ds^{2}}=0$ for some points $s \in (s_{m}, s_{2})$.
Choosing $\bar{s} \in (s_{m}, s_{2})$ so that $\frac{d^{2}x(\bar{s})}{ds^{2}} = 0$ and $\frac{d^{2}x(s)}{ds^{2}} > 0$ for $s \in (s_{m}, \bar{s})$, then we have that $\frac{dx(\bar{s})}{ds} > 0$ and $\frac{d^{3}x(\bar{s})}{ds^{3}} \leq 0$.
Using equation (\ref{2.3}),
we get $$0 \geq \frac{d^{3}x(\bar{s})}{ds^{3}} = \frac{n-1}{z(\bar{s})^{2}}(\frac{dz(\bar{s})}{ds})^{2}\frac{dx(\bar{s})}{ds} > 0,$$ it's contradictory.
Then, if there exists one point $s_{m} \in (s_{1},s_{2})$ so that $\frac{dx(s_{m})}{ds} = 0$,
we have that $\frac{d^{2}x(s)}{ds^{2}} > 0$ for $s \in (x_{1},x_{2})$.
\end{proof}

\par\bigskip \noindent{\bf Lemma 4.2.} {\ Suppose there exists $s_{m} \in (s_{1}, s_{2})$ so that $\frac{dx(s_{m})}{ds} = 0$, if $s \in (s_{m}, s_{2})$ and $z(s) \geq \frac{-\lambda+\sqrt{\lambda^{2}+4(n-1)}}{2}$, then $x(s) < \lambda$; if $s \in (s_{m}, s_{2})$ and $z(s) \geq \frac{\lambda+\sqrt{\lambda^{2}+4(n-1)}}{2}$, then $x(s) < 0$
.}\\
\begin{proof}
Since there exists $s_{m} \in (s_{1}, s_{2})$ so that $\frac{dx(s_{m})}{ds} = 0$.
By Lemma 4.1, we have $\frac{d^{2}x(s)}{ds^{2}} > 0$ for $s \in (s_{1}, s_{2})$ and $\frac{dx(s)}{ds} > 0$ for $s \in (s_{m}, s_{2})$.
If the second branch of the curve $\Gamma(s)$ may be written as $(\beta(z), z)$,
then $\beta''(z) = \frac{\frac{d^{2}x(s)}{ds^{2}}}{(\frac{dz(s)}{ds})^{4}}> 0$ for $z \in (z(s_{2}), z(s_{1}))$ and $\beta'(z) = \frac{\frac{dx(s)}{ds}}{\frac{dz(s)}{ds}}< 0$ for $z \in (z(s_{2}), z(s_{m}))$.
It follows from (\ref{3sec1}) that for $z \in (z(s_{2}), z(s_{m}))$,
\begin{equation}\label{3sec5}\aligned
0 < \frac{\beta''(z)}{1+\beta'(z)^2} = & \left(z - \frac{n-1}{z} \right) \beta'(z) - \beta(z) + \lambda\sqrt{1+\beta'(z)^2}\\
\leq & \left(z - \frac{n-1}{z} \right) \beta'(z) - \beta(z) + \lambda(1+\beta'(z))\\
= & \left(z - \frac{n-1}{z} +\lambda\right) \beta'(z) - (\beta(z) - \lambda),
\endaligned\end{equation}
where $\beta'(z) < 0$.

\noindent
When $z \geq \frac{-\lambda+\sqrt{\lambda^{2}+4(n-1)}}{2}$, we have $\beta(z) < \lambda$.
For $z \in (z(s_{2}), z(s_{m}))$, according to the equation (\ref{3sec1})again, we have that
\begin{equation}\label{3sec6}\aligned
0 < \frac{\beta''(z)}{1+(\beta'(z))^2} = & \left(z - \frac{n-1}{z} \right) \beta'(z) - \beta(z) + \lambda\sqrt{1+\beta'(z)^2}\\
\leq & \left(z - \frac{n-1}{z} \right) \beta'(z) - \beta(z) - \lambda\beta'(z)\\
= & \left(z - \frac{n-1}{z}-\lambda\right) \beta'(z) - \beta(z),
\endaligned\end{equation}
where $\beta'(z) < 0$.
When $z \geq \frac{\lambda+\sqrt{\lambda^{2}+4(n-1)}}{2}$, we have $\beta(z) < 0$.
\end{proof}

\par\bigskip \noindent{\bf Lemma 4.3.} {\ Suppose there exists a point $s_{m} \in (s_{1}, s_{2})$ so that $\frac{dx(s_{m})}{ds} = 0$, then $z(s_{2}) < \frac{\lambda+\sqrt{\lambda^{2}+4(n-1)}}{2}$.}\\
\begin{proof}
There exists $s_{m} \in (s_{1}, s_{2})$ so that $\frac{dx(s_{m})}{ds} = 0$.
By Lemma 4.1, we have $\frac{d^{2}x(s)}{ds^{2}} > 0$ for $s \in (s_{1}, s_{2})$ and $\frac{dx(s)}{ds} > 0$ for $s \in (s_{m}, s_{2})$.
By Lemma 4.2, we have that $x(s) <0$ when $s \in (s_{m}, s_{2})$ and $z(s) \geq \frac{\lambda+\sqrt{\lambda^{2}+4(n-1)}}{2}$.
Letting $M = -x(s_{m})$, then we have $x(s) \geq - M$ for $s \in (s_{1}, s_{2})$.
If $z(s_{m}) \leq \frac{\lambda+\sqrt{\lambda^{2}+4(n-1)}}{2}$, then $z(s_{2}) < z(s_{m}) \leq \frac{\lambda+\sqrt{\lambda^{2}+4(n-1)}}{2}$.
If $z(s_{m}) > \frac{\lambda+\sqrt{\lambda^{2}+4(n-1)}}{2}$, the second branch of the curve $\Gamma(s)$ may be written as $(\beta(z), z)$,
for any small $\varepsilon > 0$, there exists a constant $m_{\varepsilon} >0$ so that $z -\frac{n-1}{z} - \lambda > m_\varepsilon$
for $z\in [\frac{\lambda+\sqrt{\lambda^{2}+4(n-1)}}{2} + \varepsilon, z(s_{m}))$.
Using inequality (\ref{3sec6}),
we get
$$\beta'(z)  >  \frac{\beta(z)}{z - \frac{n-1}{z} - \lambda} \geq -\frac{M}{m_{\varepsilon}},$$
for $z\in [\frac{\lambda+\sqrt{\lambda^{2}+4(n-1)}}{2} + \varepsilon, z(s_{m}))$.
Therefore, $|\beta(z)|$ and $|\beta'(z)|$ are uniformly bounded for $z\in [\frac{\lambda+\sqrt{\lambda^{2}+4(n-1)}}{2} + \varepsilon, z(s_{m}))$.
By the existence theory for the differential equations, we have $z(s_{2}) < \frac{\lambda+\sqrt{\lambda^{2}+4(n-1)}}{2} + \varepsilon$.
In fact, if $z(s_{2}) \geq \frac{\lambda+\sqrt{\lambda^{2}+4(n-1)}}{2} + \varepsilon$,
then $\lim_{z \to z(s_{2})} \beta'(z) = -\infty$ is contradictory with uniform boundedness of $|\beta'(z)|$ on $[\frac{\lambda+\sqrt{\lambda^{2}+2(n-1)}}{2} + \varepsilon, z(s_{m}))$.
Taking $\varepsilon \to 0$, we have $$z(s_{2}) \leq \frac{\lambda+\sqrt{\lambda^{2}+4(n-1)}}{2}.$$
Next, we will show that $z(s_{2}) < \frac{\lambda+\sqrt{\lambda^{2}+4(n-1)}}{2}$.
Supposing $z(s_{2}) = \frac{\lambda+\sqrt{\lambda^{2}+4(n-1)}}{2}$, then $x(s_{2}) \in [-M, 0]$.
Since $\frac{dz(s_{2})}{ds} = 0$ and $\frac{dx(s_{2})}{ds} = 1$,
by the existence and uniqueness of solutions for the differential equations, we deduce that $z(s)$ is the constant function $z(s) \equiv \frac{\lambda+\sqrt{\lambda^{2}+4(n-1)}}{2}$, it's contradictory.
\end{proof}

\par\bigskip \noindent{\bf Lemma 4.4.} {\ $\lim_{s \to s_{2}}x(s) < \infty$.}\\
\begin{proof}
Supposing $\lim_{s \to s_{2}}x(s) = \infty$. Since $\frac{dx(s)}{ds} < 0$ for $s$ close to $s_{1}$, there exists a point $s_{m} \in (s_{1}, s_{2})$ so that $\frac{dx(s_{m})}{ds} = 0$. By Lemma 4.1, we have that $\frac{d^{2}x(s)}{ds^{2}} > 0$ for $s \in (s_{1}, s_{2})$ and $x(s_{m}) < \lambda$.
It follows that there exists a point $s_{0} \in (s_{m}, s_{2})$ so that $x(s_{0}) = 0$.
By Lemma 4.2, we have $z(s_{0}) < \frac{\lambda+\sqrt{\lambda^{2}+4(n-1)}}{2}$.
By Lemma 4.3, we have $z(s_{2}) < \frac{\lambda+\sqrt{\lambda^{2}+4(n-1)}}{2}$. In particular, when $z \in (z(s_{2}), z(s_{0}))$, there exists some $M>0$ so that $z(s) - \frac{n-1}{z(s)} - \lambda  \leq - \frac{1}{2M}$ .

\noindent
When the second branch of the curve $\Gamma(s)$ may be written as $(\beta(z), z)$, there exists $k > 0$ so that $\beta'(z(s_{0})) = -k$.
Letting $f=-\beta'$, we have that $f(z) > 0$ for $z \in (z(s_{2}), z(s_{0}))$ and $f'(z)<0$. Using equation (\ref{3sec2}), we get
\begin{equation}\aligned
f''(z) = & \frac{2f(z)f'(z)^{2}}{1+f(z)^{2}}+\left(z - \frac{n-1}{z} \right)f'(z)(1+f(z)^{2}) \\
& + \frac{(n-1)f(z)}{z^{2}}(1+f(z)^{2}) - \lambda f'(z)f(z)\sqrt{1+f(z)^{2}} \\
\geq & \left(z - \frac{n-1}{z} \right)f'(z)(1+f(z)^{2}) - \lambda f'(z)f(z)\sqrt{1+f(z)^{2}} \nonumber \\
\geq & \left(z - \frac{n-1}{z} - \lambda \right)f'(x)(1+f(z)^{2}) \\
\geq & \left(z - \frac{n-1}{z} - \lambda \right)f'(x)f(z)^{2} \\
\geq & -\frac{1}{2M}f'(z)f(z)^{2},
\endaligned\end{equation}
Where $z \in (z(s_{2}), z(s_{0}))$. For small $\varepsilon>0$, letting $$g_{\varepsilon}(z)= \frac{k\sqrt{z(s_{0}) - (z(s_{2}) + \varepsilon) } + \sqrt{3M} }{\sqrt{z - (z(s_{2}) + \varepsilon) }} ,$$ for $z \in (z(s_{2}) + \varepsilon, z(s_{0}))$, then we get
$$g_{\varepsilon}' = -\frac{1}{2}\frac{k\sqrt{z(s_{0}) - (z(s_{2}) + \varepsilon) } + \sqrt{3M}}{(z - (z(s_{2}) + \varepsilon))^{3/2}} < 0$$
and
$$g_{\varepsilon}'' = -\frac{3}{2} \frac{1}{(k\sqrt{z(s_{0}) - (z(s_{2}) + \varepsilon) } + \sqrt{3M})^2}g_{\varepsilon}'\cdot g_{\varepsilon}^2 \leq -\frac{1}{2M}g_{\varepsilon}'\cdot g_{\varepsilon}^2.$$

\noindent
Next, the purpose is to prove $f \leq g_\varepsilon$. It is known that

$$k+\frac{\sqrt{3M} }{\sqrt{z(s_{0}) - (z(s_{2}) + \varepsilon) }}=g_\varepsilon(z(s_{0})) > f(z(s_{0}))=k$$ and $$\infty=\lim_{z \to z(s_{2})+\varepsilon}g_{\varepsilon}(z)> f(z(s_{2})+\varepsilon).$$

\noindent
Therefore, if there are some points on $(z(s_{2})+\varepsilon, z(s_{0}))$ so that $f > g_{\varepsilon}$, then $f - g_{\varepsilon}$ achieves a positive maximum at point $\tilde{z} \in (z(s_{2})+\varepsilon, z(s_{0}))$.
This leads to $(f-g_{\varepsilon})'(\tilde{z}) = 0$ and $(f-g_{\varepsilon})''(\tilde{z}) \leq 0$. Then we have $$0 \geq (f-g_\varepsilon)''(\tilde{z}) \geq  -\frac{1}{2M}f'(\tilde{z}) \left( f(\tilde{z})^2-g_\varepsilon(\tilde{z})^2 \right) > 0,$$ it's contradictory.
Therefore, we have $f \leq g_{\varepsilon}$. Taking $\varepsilon \to 0$, we have the estimate $$f(z) \leq \frac{ k\sqrt{z(s_{0}) - z(s_{2})} + \sqrt{3M} }{\sqrt{z - z(s_{2})}} ,$$ for $z \in (z(s_{2}), z(s_{0}))$.
Integrating this inequality from $z$ to $z(s_{0})$, $$\beta(z) - \beta(z(s_{0})) \leq 2 \left( k\sqrt{z(s_{0}) - z(s_{2})} + \sqrt{3M}  \right) \left( \sqrt{z(s_{0})-z(s_{2})} - \sqrt{z - z(s_{2})} \right).
$$ Since $\beta(z(s_{0}))=0$, it follows that
$$\lim_{z \to z(s_{2})}\beta(z) \leq 2 \left(\sqrt{z(s_{0})-z(s_{2})} \right)\left(k\sqrt{z(s_{0}) - z(s_{2})} + \sqrt{3M}  \right).$$
Then we have $\lim_{z \to z(s_{2})} \beta(z) < \infty$
\end{proof}

\noindent
According to Lemma 4.4, there exists a point $s_{m} \in (s_{1}, s_{2})$ so that $\frac{dx(s_{m})}{ds} = 0$, then $z(s_{2}) > 0$. At this time, $\frac{dz(s_{2})}{ds} = 0$ and $\frac{dx(s_{2})}{ds} = 1$.

\par\bigskip \noindent{\bf Lemma 4.5.} {\ Suppose there exists a point $s_{m} \in (s_{1}, s_{2})$ so that $\frac{dx(s_{m})}{ds} = 0$, then $z(s_{2}) > 0$.}\\
\begin{proof}
If there exists a point $s_{m} \in (s_{1}, s_{2})$ so that $\frac{dx(s_{m})}{ds} = 0$, then we have $\frac{d^{2}x(s)}{ds^{2}} > 0$.
When the second branch of the curve $\Gamma(s)$ may be written as $(\beta(z), z)$,
we have $ -1 < \frac{\beta'(z)}{\sqrt{1 + \beta'(z)^2}} < 0$  for $z \in (z(s_{2}), z(s_{m}))$, and then there exists $\epsilon>0$ and $z_{\epsilon} \in (z(s_{2}), z(s_{m}))$ so that $\frac{\beta'(z_{\epsilon})}{\sqrt{1 + \beta'(z_{\epsilon})^2}} = -\epsilon$. Besides, by Lemma 4.4, there exists $M \geq 0$ so that $\beta(z) < M$ for $z \in (z(s_{2}), z(s_{1}))$.

\noindent
Letting $\theta(z) = \arctan \beta'(z)$, for $z \in (z(s_{2}), z_{\epsilon})$, we get that
\begin{equation}\aligned
\frac{d}{dz} \left( \ln \sin \theta(z) \right)
= & \frac{1}{\beta'(z)}\frac{\beta''(z)}{1+\beta'(z)^2} \\
= & \left(z - \frac{n-1}{z} \right) - \frac{\beta(z)}{\beta'(z)} + \frac{\lambda\sqrt{1+\beta'(z)^2}}{\beta'(z)}\\
\leq & \left(z - \frac{n-1}{z}\right)  - \frac{M}{\beta'(z_{\epsilon})} + \frac{\lambda}{\beta'(z)}(1-\beta'(z))\\
\leq & \left(z - \frac{n-1}{z} - \lambda\right)  - \frac{M-\lambda}{\beta'(z_{\epsilon})},\\ \nonumber
\endaligned\end{equation}
where $\beta'(z)<0$.
Integrating this inequality from $z$ to $z_{\epsilon}$, we get
\begin{equation}\aligned\label{3sec7}
\ln \left(\frac{\sin \theta(z_{\epsilon})}{\sin \theta(z)} \right) \leq \frac{1}{2}(z_{\epsilon})^2 + (n-1)\ln \left(\frac{z}{z_{\epsilon}} \right)
- \lambda z_{\epsilon}+ \frac{(M-\lambda) z_{\epsilon}}{(-\beta'(z_{\epsilon}))}.
\endaligned\end{equation}
Since $\sin \theta(z_{\epsilon}) = \frac{\beta'(z_{\epsilon})}{\sqrt{1 + \beta'(z_{\epsilon})^2}} = -\epsilon$, taking $z \to z(s_{2})$ for (\ref{3sec7}), we get
$$0 < \epsilon \leq (-\sin \theta(z(s_{2})))\cdot\left( \frac{z(s_{2})}{z_{\epsilon}} \right)^{n-1} e^{\frac{1}{2}(z_{\epsilon})^2 - \lambda z_{\epsilon}+\frac{(M-\lambda) z_{\epsilon}}{(-\beta'(z_{\epsilon}))}}.$$ Then $z(s_{2}) > 0$.
\end{proof}

\subsection{Behavior of the second branch of the curve $\Gamma(s)$ for small $b > -(4n+1)\lambda$}
From the above lemmas, we can get the basic description of the second branch of the curve $\Gamma(s)$: If there exists a point $s_{m} \in (s_{1}, s_{2})$ so that $\frac{dx(s_{m})}{ds} = 0$, then we have $\frac{dx^{2}(s)}{ds^{2}} > 0$ on $(s_{1}, s_{2})$ and $0 <z(s_{2})< \frac{\lambda+\sqrt{\lambda^{2}+4(n-1)}}{2}$. Next, we will study the dependence of the second branch of the curve $\Gamma(s)$ and $(x(s_{2}), z(s_{2}))$ on the initial height.

\noindent
Fixing $b \in (0, \bar{b}]$, where $\bar{b}$ is defined in Proposition 3.4.
For $b>-(4n+1)\lambda$, let $\gamma_{b}$ denote the solution of (\ref{2sec1}) with $\gamma_{b}(0) = b$ and $\gamma_{b}'(0) = 0$. Let $z(s^b_{1})$ denote the point where $\gamma_{b}$ blow-up and $x(s^b_{1}) = \gamma_{b}(z(s^b_{1}))$.
Let $\beta_{b}$ is the unique solution to (\ref{3sec1}) with $\beta(z(s^b_{1})) = x(s^b_{1})$ and $\lim _{z\to z(s^b_{1})} \beta'(z)= \infty$.
There exists $s^b_{2} > s^b_{1}$ so that $\beta_{b}$ is defined on $(z(s^b_{2}), z(s^b_{1}))$ and either blow-up as $z \to z(s^b_{2})$ or $z(s^b_{2}) = 0$.

\par\bigskip \noindent{\bf Lemma 4.6.} {\ Supposing $z(s_{1}) \geq \frac{-\lambda + 2\sqrt{2} + 2\sqrt{\lambda^{2}+4n}}{2}$, then there exists
$s_{m} \in (s_{1}, s_{2})$ so that $z(s_{m}) \in [z(s_{1})-\sqrt{2}, z(s_{1}))$ and $\beta'(z(s_{m})) = 0$.}\\

\begin{proof}
Supposing $\beta'(z) > 0$ for $z \in [z(s_{1})-\sqrt{2}, z(s_{1}))$.
By Lemma 3.8, for $s\in (0, s_{1})$, if $z(s) \geq \sqrt{\frac{2n(b-\lambda)}{b+\lambda}}$, we have $x(s)= \beta(z(s))< \lambda$.
Since $b > -(4n+1)\lambda$, then $\sqrt{\frac{2n(b-\lambda)}{b+\lambda}} < 2\sqrt{n}$.
When $z(s_{1}) \geq \frac{-\lambda + 2\sqrt{2} + 2\sqrt{\lambda^{2}+4n}}{2}$, we have that $z(s_{1})-\sqrt{2} \geq \frac{-\lambda + 2\sqrt{\lambda^{2}+4n}}{2}>\frac{-\lambda + \sqrt{\lambda^{2}+4(n-1)}}{2}$ and $\beta(z(s_{1})) < \lambda$.
Using the equation (\ref{3sec1}), for $z \in [z(s_{1})-\sqrt{2}, z(s_{1}))$, we have
\begin{equation}\aligned
\frac{\beta''(z)}{1+\beta'(z)^2} = & \left(z - \frac{n-1}{z} \right) \beta'(z) - \beta(z) + \lambda\sqrt{1+\beta'(z)^2}\\
\geq & \left(z - \frac{n-1}{z} \right) \beta'(z) - \beta(z)+ \lambda(1+\beta'(z))\\
\geq & \left(z - \frac{n-1}{z} + \lambda \right) \beta'(z) - (\beta(z)-\lambda), \nonumber
\endaligned\end{equation}
then $\beta''(z) > 0$ on $[z(s_{1})-\sqrt{2}, z(s_{1}))$ and $\beta''(z(s_{1})-\sqrt{2}) \geq -(\beta(z(s_{1})-\sqrt{2})-\lambda)$.
Therefore, for $z \in [z(s_{1})-\sqrt{2}, z(s_{1}))$, using the equation (\ref{3sec2}), we have
\begin{equation}\aligned
\frac{\beta'''(z)}{1+\beta'(z)^2} = & \frac{2\beta'(z)(\beta''(z))^2}{( 1+\beta'(z)^2 )^2}+\left(z
-\frac{n-1}{z} \right)\beta''(z)+\frac{n-1}{z^2} \beta'(z)+\frac{\lambda\beta'(z)\beta''(z)}{\sqrt{1+\beta'(z)^2}}\\
\geq & (z - \frac{n-1}{z}) \beta''(z)+\frac{\lambda\beta'(z)\beta''(z)}{\sqrt{1+\beta'(z)^2}}\\
\geq & (z - \frac{n-1}{z} + \lambda) \beta''(z), \nonumber
\endaligned\end{equation}
then $\beta'''(z) > 0$ on $[z(s_{1})-\sqrt{2}, z(s_{1}))$. For $z \in [z(s_{1})-\sqrt{2}, z(s_{1}))$, integrating this inequality repeatedly from $z(s_{1})-\sqrt{2}$ to $z$, we get
$$\beta(z) \geq \beta(z(s_{1})-\sqrt{2}) - \frac{1}{2}(\beta(z(s_{1})-\sqrt{2}) - \lambda) \left(z- (z(s_{1})-\sqrt{2}) \right)^2,$$ where we use $\beta''(z(s_{1})-\sqrt{2}) \geq -(\beta(z(s_{1})-\sqrt{2})-\lambda)$. This shows us that $\beta(z(s_{1})) \geq \lambda$, it's contradictory.
\end{proof}

\noindent
Let $\bar{b} > 0$ be the value given in the conclusion of Proposition 3.4 so that if $b \in (0,\bar{b}]$, $\gamma(z) \leq 0$ for $z\in[\sqrt{\frac{2nb}{b+\lambda}}, z(s_{1}))$. We also assume that $\bar{b}$ is chosen so small that $z(s_{m}) > 2\sqrt{2n}$ and $x(s_{1}) \geq -\frac{1}{8\sqrt{2n}} > -\frac{7n+1}{4\sqrt{2n}}$ when $b \in (0, \bar{b}]$. So we will get the following conclusion.

\par\bigskip \noindent{\bf Lemma 4.7.} {\ If $b \in (0,\bar{b}]$, then $2x(s_{1}) \leq \beta(z) < 0$ for $z \in [2\sqrt{2n}, z(s_{1})]$.}\\
\begin{proof}
Let $\alpha(x)$ denote the curve that connects $\gamma$ and $\beta$ near the point $(x(s_{1}), z(s_{1}))$, then $\alpha$ is a solution of (\ref{3sec3}) with $\alpha(x(s_{1})) = z(s_{1})$ and $\alpha'(x(s_{1})) = 0$.
Using equation (\ref{3sec4}), we have $\alpha'''(x(s_{1})) = x(s_{1})\alpha''(x(s_{1}))$. Therefore, we get Taylor expansion of $\alpha(x)$,
\begin{equation}\aligned\label{3sec8}
\alpha(x) = & z(s_{1}) + \frac{1}{2}\alpha''(x(s_{1}))(x-x(s_{1}))^2\\
 & + \frac{1}{6}x(s_{1})\alpha''(x(s_{1}))(x-x(s_{1}))^3 + \mathcal{O}(|x-x(s_{1})|^4),
\endaligned\end{equation} as $x \to x(s_{1})$.

\noindent
For $z<z(s_{1})$ and near $z(s_{1})$, the curve $\Gamma$ is concave down and heads back towards the $x$-axis.
We can find $s,\ \ t >0$ so that $$\alpha(x(s_{1}) + t)= z = \alpha(x(s_{1}) - s).$$
Using (\ref{3sec8}) and $x(s_{1})<0$, we can get $t>s$.

\noindent
To prove the lemma, we consider the function $\delta(z) = \gamma(z) + \beta(z)$. Since $\gamma(z) \leq 0$ and $\beta(z) < 0$ for $z \in [2\sqrt{2n}, z(s_{1}))$, we have $\delta(z) < 0$. By the previous discussion, then $$\delta(z) = (x(s_{1})+ t) + (x(s_{1}) -s)  > 2x(s_{1}) = \delta(z(s_{1})),$$ for $z < z(s_{1})$ and near $z(s_{1})$.

\noindent
Now, we will prove that $\delta(z) > 2x(s_{1})$ on $[2\sqrt{2n}, z(s_{1}))$. Supposing $\delta(z) = 2x(s_{1})$ for some $z \in [2\sqrt{2n}, z(s_{1}))$.
Choosing $\tilde{z} \in [2\sqrt{2n}, z(s_{1}))$ so that $\delta(\tilde{z}) = 2x(s_{1})$.
It follows from the previous discussion that there exists a point $\bar{z} \in (\tilde{z}, z(s_{1}))$ so that $\delta$ achieves a negative maximum. This leads to $\delta'(\bar{z}) = 0$ and $\delta''(\bar{z}) \leq 0$. Therefore, we get
$$0 \geq \frac{\delta''(\bar{z})}{1 + \gamma'(\bar{z})^2} = -\delta(\bar{z}) > 0,$$ it's contradictory.
So we have that $\delta > 2x(s_{1})$ on $[2\sqrt{2n}, z(s_{1}))$. Since $\gamma(z) \leq 0$ on $[2\sqrt{2n}, z(s_{1}))$, we have that $ 2x(s_{1}) \leq \gamma(z) + \beta(z) \leq \beta(z) < 0$ on $[2\sqrt{2n}, z(s_{1}))$.
\end{proof}

\par\bigskip \noindent{\bf Lemma 4.8.} {\ Letting $b \in (0, \bar{b}]$. If $\beta<0$ on $(z(s_{2}), z(s_{1}))$, then
$$z(s_{2}) \leq \frac{8(n-1)}{(\pi-2)+\frac{16n-1}{\sqrt{2n}}(-\lambda)}(-x(s_{1})).$$}\\
\begin{proof}
According to our previous assumptions on $\bar{b}$, we have that $z(s_{m})> 2\sqrt{2n}$, $\beta'' > 0$ and $\beta(2\sqrt{2n}) \geq 2x(s_{1})$.
By Lemma 4.2, we have $2x(s_{1}) \leq \beta(2\sqrt{2n}) < \lambda$.
For $z \in (z(s_{2}), z(s_{m}))$, using the inequality (\ref{3sec5}) and choosing$z=2\sqrt{2n}$, we have
$$\beta'(2\sqrt{2n}) > \frac{\beta(2\sqrt{2n})-\lambda}{2\sqrt{2n} - \frac{n-1}{2\sqrt{2n}} + \lambda}
\geq \frac{2x(s_{1})-\lambda}{\frac{7n+1}{2\sqrt{2n}} + \lambda} > -1.$$
For $z \in (z(s_{2}), 2\sqrt{2n})$,
we have that
\begin{equation}\aligned
\frac{d}{dz} \left(\arctan\beta'(z) \right) = & \left(z - \frac{n-1}{z} \right) \beta'(z) - \beta(z) + \lambda\sqrt{1+\beta'(z)^2}\\
\leq & - \frac{n-1}{z}\beta'(z) - \beta(z) + \lambda(1+\beta'(z)) \\
\leq & (- \frac{n-1}{z(s_{2})}+\lambda)\beta'(z)+(\lambda-\beta(2\sqrt{2n})), \nonumber
\endaligned\end{equation}
where $\beta'(z)<0$.
Integrating this inequality from $z(s_{2})$ to $2\sqrt{2n}$, then
\begin{equation}\aligned
\frac{\pi}{4}
\leq & \arctan\beta'(2\sqrt{2n})+\frac{\pi}{2} \\
\leq & (-\frac{n-1}{x(s_{2})}+\lambda)\beta(2\sqrt{2n})+2\sqrt{2n}(\lambda-\beta(2\sqrt{2n}))\\
= & (-\beta(2\sqrt{2n}))(\frac{n-1}{z(s_{2})}-\lambda+2\sqrt{2n})+2\sqrt{2n}\lambda\\
\leq & (-2x(s_{1}))(\frac{n-1}{z(s_{2})}-\lambda+2\sqrt{2n})+2\sqrt{2n}\lambda, \nonumber
\endaligned\end{equation}
where $\beta(z(s_{2}))<0$, $2x(s_{1}) \leq \beta(2\sqrt{2n}) < \lambda$ and $\beta'(2\sqrt{2n}) > -1$.
Rearranging this inequality to estimate $z(s_{2})$ and using $x(s_{1}) \geq -\frac{1}{8\sqrt{2n}}$,
we get $$z(s_{2}) \leq \frac{8(n-1)}{(\pi-2)+\frac{16n-1}{\sqrt{2n}}(-\lambda)}(-x(s_{1})).$$
\end{proof}

\par\bigskip \noindent{\bf Lemma 4.9.} {\ There exists $\bar b>0$ so that for $b \in (0, \bar b]$, there is a point $z(s^b_{m}) \in (z(s^b_{2}), z(s^b_{1}))$ so that $\beta_{b}'(z(s^b_{m})) =0$ and $0 < \beta_{b}(z(s^b_{2})) < \infty$.}\\
\begin{proof}
Fixing $b \in (0, \bar{b}]$ and leting $\beta = \beta_{b}$. By Lemma 4.6, we have that there exists $s_{m} \in (s_{1}, s_{2})$ so that $z(s_{m}) > 2\sqrt{2n}$ and $\beta'(z(s_{m})) = 0$. Therefore, we have $\beta''>0$ on $(z(s_{2}), z(s_{1}))$ and $\beta' < 0$ on $(z(s_{2}), z(s_{m}))$. Besides, by Lemma 4.3, we get $\beta(z(s_{2})) < \infty$.

\noindent
Supposing $\beta(z(s_{2})) \leq 0$ so that $\beta<0$ on $(z(s_{2}), z(s_{m}))$.
By Lemma 4.8, there exists $0<\epsilon<\frac{\lambda + \sqrt{\lambda^{2} + (n-1)}}{2}$ so that $\beta'(\epsilon)=-1$ for sufficiently small $b$. Therefore, using equation (\ref{3sec2}), we get
\begin{equation}\aligned
\frac{\beta'''(z)}{1+\beta'(z)^2} = & \frac{2\beta'(z)(\beta''(z))^2}{( 1+\beta'(z)^2 )^2}+\left(z -\frac{n-1}{z} \right)\beta''(z)+\frac{n-1}{z^2} \beta'(z)\\
&+\frac{\lambda\beta'(z)\beta''(z)}{\sqrt{1+\beta'(z)^2}}\\
\leq & (z - \frac{n-1}{z}-\lambda) \beta''(z), \nonumber
\endaligned\end{equation}
then $\beta'''(z)<0$ for $z \in (z(s_{2}), \frac{\lambda+\sqrt{\lambda^{2}+4(n-1)}}{2}]$.
Using the equation (\ref{3sec1}), then
\begin{equation}\aligned
\beta''(z)
\geq & \frac{\beta''(z)}{1+\beta'(z)^2} \\
= & \left(z - \frac{n-1}{z} \right) \beta'(z) - \beta(z) + \lambda\sqrt{1+\beta'(z)^2}\\
\geq & \left(z - \frac{n-1}{z} \right) \beta'(z) - \beta(z) + \lambda(1-\beta'(z))\\
= & \left(z - \frac{n-1}{z} - \lambda\right) \beta'(z) - (\beta(z) - \lambda),\nonumber
\endaligned\end{equation}
Taking $z = \frac{\lambda + \sqrt{\lambda^{2} + (n-1)}}{2}$,
we have $\beta''(\frac{\lambda + \sqrt{\lambda^{2} + (n-1)}}{2}) \geq -(\beta(\frac{\lambda + \sqrt{\lambda^{2} + (n-1)}}{2}) - \lambda)$, and then $\beta''(z) \geq -(\beta(\frac{\lambda + \sqrt{\lambda^{2} + (n-1)}}{2}) - \lambda)$ for $z \in (z(s_{2}), \frac{\lambda + \sqrt{\lambda^{2} + (n-1)}}{2}]$. Integrating this inequality from $\epsilon$ to $\frac{\lambda + \sqrt{\lambda^{2} + (n-1)}}{2}$. In addition, $\beta'(z)<0$ for $z \in (\epsilon, \frac{\lambda + \sqrt{\lambda^{2} + (n-1)}}{2})$, we have
$$\beta'(\epsilon) \leq (\frac{\lambda + \sqrt{\lambda^{2} + (n-1)}}{2} - \epsilon) (\beta(\epsilon) - \lambda).$$
Therefore,
$$\beta(\epsilon)\geq -\frac{1}{\frac{\lambda + \sqrt{\lambda^{2} + (n-1)}}{2} - \epsilon}+\lambda.$$
Choosing $z \in (z(s_{2}), \epsilon)$, we have
\begin{equation}\aligned\label{3sec9}
\frac{\beta''(z)}{\beta'(z)} = & \left(z - \frac{n-1}{z} \right)(1+\beta'(z)^2) - \frac{\beta(z)}{\beta'(z)}(1+\beta'(z)^2)+ \frac{\lambda(1+\beta'(z)^2)^{\frac{3}{2}}}{\beta'(z)}\\
\leq & \left(z - \frac{n-1}{z} +\frac{\lambda}{\beta'(z)}(1+\beta'(z)^2)^{\frac{1}{2}}\right)(1+\beta'(z)^2)\\
\leq & \left(z - \frac{n-1}{z}-2\lambda\right)(1+\beta'(z)^2)\\
\leq &z - \frac{n-1}{z}-2\lambda,
\endaligned\end{equation}
where $\beta'(z) < \beta'(\epsilon) = -1$.
Integrating (\ref{3sec9}) repeatedly from $z$ to $\epsilon$, we get
$$\beta(z) \geq \beta(\epsilon) - \beta'(\epsilon) e^{-(\frac{\epsilon^{2}}{2}-2\lambda\epsilon)}\int^{\epsilon}_{z} (\frac{\epsilon}{t})^{n-1}dt= \beta(\epsilon)+e^{-(\frac{\epsilon^{2}}{2}-2\lambda\epsilon)}\int^{\epsilon}_{z} (\frac{\epsilon}{t})^{n-1}dt.$$
Combining this with $\beta(\epsilon)\geq -\frac{1}{\frac{\lambda + \sqrt{\lambda^{2} + (n-1)}}{2} - \epsilon}+\lambda$, we have that
$$0 > \beta(z) \geq e^{-(\frac{\epsilon^{2}}{2}-2\lambda\epsilon)}\int^{\epsilon}_{z} (\frac{\epsilon}{t})^{n-1}dt-\frac{1}{\frac{\lambda+\sqrt{\lambda^{2}+(n-1)}}{2}-\epsilon}+\lambda.$$
By Proposition 3.4 and Lemma 4.8, we have that $z(s_{2})$ is sufficiently close to $0$ when $b > -(4n+1)\lambda$ is sufficiently small, and then $\beta(z) \geq 0$ for $z\to z(s_{2})$, it's contradictory.
Therefore, $0 < \beta_b(z(s^b_{2})) < \infty$ for sufficiently small $b>-(4n+1)\lambda$.
\end{proof}

\vspace*{2mm}
%%%%%%%%%%%%%%%%%%%%%%%%%%%%%%%%%
\section{An Immersed $S^{n}$ $\lambda$-hypersurface}

In this section, we will complete the proof of Theorem 1.1. We consider the set
$$
T =  \{\tilde{b} \, : \, \forall b \in (0, \tilde{b}], \exists s_{m}^b \in (s^b_{2}, s^b_{1}) \textrm{ so that }
\frac{dx(s^b_{m})}{ds} = 0 \textrm{ and } x(s^b_{2}) > 0\}.
$$
By Lemma 4.9, we have that this set is not empty. Following Angenent's argument in~\cite{A}, let $b_{0}$ be the supremum of this set: $$b_{0} = \sup_{\tilde{b}} T.$$
In~\cite{CW}, we know that the curve $\Gamma(s)=(b\cos\frac{s}{b}, b\sin\frac{s}{b})$ is a circle, where $b=\frac{-\lambda+\sqrt{\lambda^{2}+4n}}{2}$, which is special solution of (\ref{2.2}). It's obvious that $b_{0} \leq \frac{-\lambda+\sqrt{\lambda^{2}+4n}}{2}$. We want to show $\Gamma_{b_{0}}$ intersects the $x$-axis perpendicularly at $z(s^{b_{0}}_{2})$ ($x(s^{b_{0}}_{2}) = 0$).

\par\bigskip \noindent{\bf Lemma 5.1.} {\ $b_{0}< \frac{-\lambda+\sqrt{\lambda^{2}+4n}}{2}$.}\\
\begin{proof}
Suppose $\frac{dx(s)}{ds}<0$ and $\frac{d^{2}x(s)}{ds^{2}} \geq 0$ on $(s^{b_{0}}_{1}, s^{b_{0}}_{2})$.
It writes the second branch of the curve $\Gamma_{b_{0}}(s)$ as a form of $(\beta_{b_{0}}(z), z)$,
we have that $\beta_{b_{0}}'>0$ and $\beta''_{b_{0}} \geq 0$ on $(z(s^{b_{0}}_{2}), z(s^{b_{0}}_{1}))$.
Then, the curve $\beta_{b_{0}}$ must intersect with the $x$-axis
$(z(s^{b_{0}}_{2}) = 0)$ and there exists $m>0$ so that $\beta_{b_{0}}(z) \leq -m$ for $z \in (0,\frac{1}{2}]$.
Let $\{b_{n}\}$ be an increasing sequence that converges to $b_{0}$. Let $\beta_{n}$ be the solution of (\ref{3sec1}) corresponding to the initial height $b_{n}$. Fixing $0 < \varepsilon< \frac{1}{2}$, by continuity of the curve $\Gamma_{b}(s)$, there exists $N=N(\varepsilon)>0$ so that for $n>N$, we have that
$|z(s^{n}_{2}) - z(s^{b_{0}}_{2})| < \varepsilon$ and $|\beta_n(\varepsilon) - \beta_{b_{0}}(\varepsilon)|< \varepsilon$. Therefore, $z(s^{n}_{2}) < \varepsilon$ and $\beta_n(\varepsilon) < -m$,
it follows that $\beta_{n}$ intersects the $z$-axis at point $z(s^n_{0}) < \varepsilon$.
Besides, the curve $(\beta_{n}(z), z)$ can be written as the curve $(x,\alpha_{n}(x))$ for $x \in [-m/2,0]$,
it follows from (\ref{2.2}) that
\begin{equation}\aligned\label{4sec1}
\frac{\alpha_{n}''(x)}{1+\alpha_{n}'(x)^2} = \left(\frac{n-1}{\alpha_{n}} - \alpha_{n}\right) + x\alpha_{n}'(x) + \lambda\sqrt{1+\alpha_{n}'(x)^2},
\endaligned\end{equation}
where $\alpha_{n}''(x) = \frac{\frac{d^{2}z(s)}{ds^{2}}}{(\frac{dx(s)}{ds})^{4}}$ and $\frac{dx(s)}{ds} = \frac{1}{\sqrt{1+\alpha_{n}'(x)^{2}}}$.

\noindent
In fact, since $\beta_{n}(\varepsilon) < -m < -\frac{m}{2}$, $\alpha'_{n}(x)=\frac{1}{\beta'_{n}(z)} < 0$ and $\alpha''_{n}(x)=\frac{-\beta''_{n}(z)}{\beta'_{n}(z)^{3}}$,
we have the following estimates $$0 < \alpha_{n}(x) < \alpha_{n}(-\frac{m}{2}) < \varepsilon,\ \ \ \  -\frac{2\varepsilon}{m}\leq\alpha_{n}'(x) < 0,\ \ \ \ \alpha''_{n}(x) \geq 0, \ \ \ x \in [-m/2,0].$$
Using equation (\ref{4sec1}), for $x \in [-m/2,0]$, then
\begin{equation}\aligned
\alpha_{n}''(x)
\geq& \frac{\alpha_{n}''(x)}{1+\alpha_{n}'(x)^2}\\
= & \left(\frac{n-1}{\alpha_{n}} - \alpha_{n}\right) + x\alpha_{n}'(x) + \lambda\sqrt{1+\alpha_{n}'(x)^2}\\
\geq & \left(\frac{n-1}{\alpha_{n}} - \alpha_{n}\right) + x\alpha_{n}'(x) + \lambda(1-\alpha_{n}'(x))\\
\geq & \frac{n-1}{\alpha_{n}} - \alpha_{n} + \lambda(1+\frac{2\varepsilon}{m})\\
\geq & \frac{n-1}{\varepsilon} - \varepsilon + \lambda(1+\frac{2\varepsilon}{m}), \nonumber
\endaligned\end{equation}
where $0 < \alpha_{n}(x) < \varepsilon$ and $-\frac{2\varepsilon}{m}\leq\alpha_{n}'(x) < 0$.
For small $\varepsilon>0$, we have $$\alpha_{n}''(x) \geq \frac{n-1}{2\varepsilon} + \lambda(1+\frac{2\varepsilon}{m})$$
Integrating this inequality repeatedly from $x$ to $0$, then
$$\alpha_{n}(x) \geq \alpha_{n}(0) + \alpha_{n}'(0)x + \frac{1}{2}(\frac{n-1}{2\varepsilon} + \lambda(1+\frac{2\varepsilon}{m}))x^2 \geq \frac{1}{2}(\frac{n-1}{2\varepsilon} + \lambda(1+\frac{2\varepsilon}{m}))x^2.$$
Taking $x=-\frac{m}{4}$, we have $\alpha_n(-\frac{m}{4}) \geq \frac{m^2}{32}(\frac{n-1}{2\varepsilon} + \lambda(1+\frac{2\varepsilon}{m}))$.
Since there exists the point $\bar{z}\in (z(s^n_{2}), \varepsilon)$ so that $\beta_n(\bar{z}) = -\frac{m}{4}$, we have that $\frac{m^2}{32}(\frac{n-1}{2\varepsilon} + \lambda(1+\frac{2\varepsilon}{m})) \leq \bar{z} < \varepsilon$. This contradicts the fact that if $\varepsilon > 0$ is sufficiently small, where $\lambda$ is bounded negative.
Then, there dose not exist $\frac{dx(s)}{ds}<0$ and $\frac{d^{2}x(s)}{ds^{2}} \geq 0$ on $(s^{b_{0}}_{1}, s^{b_{0}}_{2})$.
When $b_{0} = \frac{-\lambda+\sqrt{\lambda^{2}+4n}}{2}$, we have that $x(s)=b_{0}\cos\frac{s}{b_{0}} < 0$ and $z(s)=b_{0}\sin\frac{s}{b_{0}} > 0$ for $s\in (s^{b_{0}}_{1}, s^{b_{0}}_{2})$, then $\frac{dx(s)}{ds}=-\sin\frac{s}{b_{0}} < 0$ and $\frac{d^{2}x(s)}{ds^{2}}=-\frac{1}{b_{0}}\cos\frac{s}{b_{0}} > 0 $ for $s\in (s^{b_{0}}_{1}, s^{b_{0}}_{2})$. So we get $b_{0} < \frac{-\lambda+\sqrt{\lambda^{2}+4n}}{2}$.
\end{proof}

\begin{proof}[Proof of Theorem 1.1]
The proof of the theorem 1.1 is divided into two parts. Part one, we will show that there exists $s^{b_{0}}_{m} \in (s^{b_{0}}_{1}, s^{b_{0}}_{2})$ so that $\frac{dx(s^{b_{0}}_{m})}{ds} = 0$.
Supposing that $\frac{dx(s)}{ds} < 0$ on $(s^{b_{0}}_{1}, s^{b_{0}}_{2})$. By Lemma 5.1, we have that it is impossible to have $\frac{dx(s)}{ds}<0$ and $\frac{d^{2}x(s)}{ds^{2}} \geq 0$ on $(s^{b_{0}}_{1}, s^{b_{0}}_{2})$. Therefore, there exists $\frac{d^{2}x(s)}{ds^{2}} < 0$ for near $s^{b_{0}}_{2}$.
Let $\{b_{n}\}$ be an increasing sequence that converges to $b_{0}$. Let $\Gamma_{n}(s)$ be the solution of (\ref{2.2}) corresponding to the initial height $b_{n}$. According to the definition of $b_{0}$, we have that for sufficiently large $n$, $\frac{d^{2}x(s)}{ds^{2}} > 0$ on $(s^{b_{n}}_{1}, s^{b_{n}}_{2})$,  which contradicts the fact that $\frac{d^{2}x(s)}{ds^{2}} < 0$ for near $s^{b_{0}}_{2}$.

\noindent
Part two, we will show that $x(s^{b_{0}}_{2}) = 0$.
According to the discussion of the previous paragraph, we have that $z(s^{b_{0}}_{2}) > 0$ and $|x(s^{b_{0}}_{2})| < \infty$. Therefore, by the continuous dependence of the curve $\Gamma_{b}(s)$ on the initial height, there is a $\bar{\delta}>0$ so that when $|b-b_{0}| < \bar{\delta}$ and $b > 0$  there exists $s^{b}_{m} \in (s^{b}_{1}, s^{b}_{2})$ so that $\frac{dx(s^{b}_{m})}{ds} = 0$.
To complete the proof, we will discuss three cases: $x(s^{b_{0}}_{2}) > 0$, $x(s^{b_{0}}_{2}) < 0$ and $x(s^{b_{0}}_{2}) = 0$.
Our purpose is to prove that the first and the second case cannot occur.
If $x(s^{b_{0}}_{2}) > 0$, it follows from the continuity of the curve $\Gamma_{b}(s)$ that we can find $\delta_{1} \in (0, \bar{\delta})$ so that when $|b-b_{0}|< \delta_{1}$ and $b > 0$ the curve $\Gamma_{b}$ has a vertical tangent point $(x(s^{b}_{2}), z(s^{b}_{2}))$ with $x(s^{b}_{2}) > 0$. Then $b_{0} + \delta_{1}/2 \in T$, this contradicts the definition of $b_{0}$. If $x(s^{b_{0}}_{2}) < 0$, it follows from the continuity of the curve $\Gamma_{b}(s)$ that there exists $\delta_{2} \in (0, \bar{\delta})$ so that when $|b-b_{0}|< \delta_{2}$ and $b > 0$ the curve $\Gamma_{b}$ has a vertical tangent point $(x(s^{b}_{2}), z(s^{b}_{2}))$ with $x(s^{b}_{2}) < 0$. However, $b_{0} - \delta_{2}/2 \in T$, this also contradicts the definition of $T$.

\noindent
From the above we have that the curve $\Gamma_{b_{0}}(s)$ intersects the $z$-axis perpendicularly at the point $(0, z(s^{b_{0}}_{2}))$ where $z(s^{b_{0}}_{2}) \in (0, \frac{\lambda+\sqrt{\lambda^{2}+4(n-1)}}{2})$.
Let $\Gamma_{b_{0}}$ be the solution to (\ref{2.2}) obtained by shooting perpendicularly to the $x$-axis at the point $(b_{0}, 0)$ in the right-half plane.
We show that $\Gamma_{b_{0}}$ follows the curve $(\gamma_{b_{0}}(z), z)$ from the $x$-axis, across the
$z$-axis, to the point $(x(s^{b_{0}}_{1}), z(s^{b_{0}}_{1}))$, and then it follows the curve $(\beta_{b_{0}}(z), z)$ until it intersects the $z$-axis perpendicularly at $(0,z(s^{b_{0}}_{2}))$.
Since the equation is symmetric with respect to reflections across the $z$-axis, we see that $\Gamma_{b_{0}}$ continues along the
reflected curves $(-\beta_{b_{0}}(z), z)$ and $(-\gamma_{b_{0}}(z), z)$ until it exits the right-half plane, where
it intersects the $x$-axis perpendicularly at the point $(-b_{0},0)$.
That is, $\Gamma_{b_{0}} = \gamma_{b_{0}} \cup \beta_{b_{0}} \cup -\beta_{b_{0}} \cup -\gamma_{b_{0}}$.
Since the rotation of $\Gamma_{b_{0}}$ about the $x$-axis is smooth in a neighborhood of these symmetric points and $\Gamma_{b_{0}}$ is smooth in the vertical tangent points $z(s^{b_{0}}_{1})$ where the $\gamma$ and $\beta$ curves meet, then the rotation of $\Gamma_{b_{0}}$ about the $x$-axis is a smooth manifold.
Notice the facts that $(\gamma_{b_{0}}(z), z)$ and $(-\beta_{b_{0}}(z), z)$ intersect transversally and the shapes of $\gamma$ and $\beta$ are convex, we have that the curve $\Gamma_{b_{0}}$ has two self-intersections.
Therefore, its rotation about the $x$-axis is an immersed, non-embedded $S^n$ $\lambda$-hypersurface in $\mathbb{R}^{n+1}$.
\end{proof}

\vspace*{2mm}
%%%%%%%%%%%%%%%%%%%%%%%%%%%%%%%%%
\bibliographystyle{Plain}

\end{document}